\documentclass[reqno]{amsart}
\usepackage[hmargin={25mm,25mm}]{geometry}
\usepackage[dvipdfmx]{graphicx,color}
\newtheorem{thm}{Theorem}[section]

\theoremstyle{definition}

\theoremstyle{remark}

\title[IPS function for Schr\"0dinger equation ]{Integrating the probe and singular sources methods: IV.  
\\
IPS function for the Schr\"odinger equation}

\author{Masaru \textsc{Ikehata}}
\address{
Professor Emeritus at Gunma University;
Professor Emeritus at Hiroshima University, Graduate School of Advanced Science and Engineering, 
Hiroshima University\\
Higashihiroshima, Japan
}
\email{ikehataprobe@gmail.com}

\thanks{
The author was partially supported by Grant-in-Aid for
Scientific Research (C)(No. 24K06812) of Japan  Society for
the Promotion of Science.
 }


\begin{document}

\begin{abstract}
The integrated theory of the probe and singular sources methods (IPS) is developed
for an inverse obstacle problem governed by the stationary Schr\"odinger equation in a bounded domain.
The unknown obstacles are penetrable, and their surface is modeled by a part of the support of the potential in the governing equation.
The main results concern an analytical detection method for these obstacles from the Dirichlet-to-Neumann map.
They consist of three parts: a singular sources method via the probe method using a solution with higher-order singularity for the governing equation
of the background medium; the discovery of an IPS function whose two ways of decomposition give us the indicator functions for both the probe and singular sources methods;
a completely integrated version of both methods,  which means their indicator functions coincide.
Furthermore, a result on Side B of IPS is also given, concerning the blowing-up property of a sequence calculated
from the Dirichlet-to-Neumann map.

\end{abstract}

\maketitle

\vskip.5cm
\noindent

\vskip.2cm
\section{Introduction}

The aim of this paper is to pursue  further the integrated theory of the probe and singular sources methods (IPS) recently developed by the author in \cite{IPS}, \cite{IPS2} and \cite{IPS3}.

In \cite{IPS} the author reconsidered the spirit of the Singular Sources Method \cite{P1} (see also \cite{NP}) and introduced its reformulation in terms of the Probe Method \cite{INew}
by considering a prototype inverse obstacle problem in a bounded domain.  Then he found a function with two decompositions, both of which
yield indicator functions for the methods.  Therein the function is called the IPS function and itself also plays the role of an indicator function.  
So the IPS function makes a bridge between the Probe and Singular Sources Methods.  
In \cite{IPS2}, the idea has been applied to an inverse obstacle problem governed by the Stokes system
in a bounded domain.
The author derived the Singular Sources Method for the Stokes system
from the Probe Method, which demonstrates an advantage of IPS.  
Furthermore, as a byproduct, it was affirmatively solved whether the Singular Sources Method has the Side B property of the Probe Method \cite{INew}, a property which states
the blowing up of a sequence calculated from the Dirichlet-to-Neumann map;
this had been a pending question in \cite{IPS}.
In \cite{IPS3}, an inverse obstacle problem governed by the Helmholtz equation in a bounded domain was considered.  Therein, the unknown obstacles consist of two types: sound-soft and sound-hard ones.
It was shown that IPS still works.  Besides, to handle such different types of obstacles, 
the author introduced a {\it method of a complementing function}.

However, all the obstacles considered above are {\it impenetrable}. 
It should be noted that at the beginning (1998), the Probe Method was also developed for two types of penetrable obstacles, as described in \cite{IProbe} and also \cite{IProbeElastic}.
One of these appears as a jump of the second-order coefficient of the governing equation---a special version of the Caler\'on problem \cite{Car},
which is the continuum model of Electrical Impedance Tomography.  A uniqueness theorem for this case was
established by Isakov \cite{Is} (1988) using his method of singular solutions.  The other case is a jump in the zeroth-order coefficient, i.e., the potential of the stationary 
Schr\"odinger equation, roughly speaking, 
which is related to the Calde\'ron problem by the Liouville transform.

This second type of problem also appears as the inverse medium problem in inverse scattering (cf. \cite{CK}) whose governing equation
is, for example, the inhomogeneous Helmholtz equation.
Though one can fully reconstruct the potential itself by the results of Nachman \cite{Na} (1988) and Novikov\cite{No} (1988), other methods,
the Factorization Method \cite{Ki} (1999) and the Singular Sources Method \cite{Pbook} (2001), have also been proposed.  
These methods, together with the Probe Method, now belong to a group of analytical methods for inverse obstacle problems.

In this paper, we focus on the second type of penetrable obstacle case and a result presented in \cite{IProbe}.
This problem is governed by the stationary Sch\"odinger equation in a bounded domain, where
the unknown obstacle appears as a jump in the potential, and the author developed the Probe Method.
However, needless to say, at that time (1998), we did not have the IPS view. 
We investigate the following question:
How does IPS work for this penetrable case?   
Our motivation stems from a purely mathematical interest in the methodology
for inverse obstacle problems governed by partial differential equations.

Now let us review the result of Section 4 in \cite{IProbe} in terms of the reformulated Probe Method \cite{INew}, which is now called the Side A of the Probe Method.
Let $\Omega$ be a bounded domain of $\mathbb R^3$ with connected boundary of class $C^{1,1}$.
Let $V$ be a real-valued function in $L^{\infty}(\Omega)$.  We assume that $0$ is not a Dirichlet eigenvalue of $-\Delta+V(x)$ in $\Omega$.
Then, given $f\in H^{\frac{1}{2}}(\partial\Omega)$ there exists a unique $u\in H^1(\Omega)$ such that

\begin{equation}
%$$
\left\{
\begin{array}{ll}
\displaystyle
-\Delta u+V(x)u=0, & x\in\Omega,
\\
\\
\displaystyle
u=f, & x\in\partial\Omega.
\end{array}
\right.
\tag {1.1}
%$$
\end{equation}

Let $\Lambda_V$ be the associated Dirichlet-to-Neumann map, that is written formally
$$\displaystyle
\Lambda_V:f\longmapsto\frac{\partial u}{\partial\nu}\vert_{\partial\Omega},
$$
where $u$ is the solution of (1.1) and $\nu$ denotes the unit outward normal to $\partial\Omega$.  The $\Lambda_Vf$ has the expression as the bounded linear functional on $H^{\frac{1}{2}}(\partial\Omega)$:
$$\displaystyle
<\Lambda_Vf,g>=\int_{\Omega}\nabla u\cdot\nabla v+V(x)uv\,dx,
$$
where $v\in H^1(\Omega)$ is an arbitrary function satisfying $v=g$ on $\partial\Omega$ in the sense of trace (cf. \cite{Gr}).

$\quad$

Throughout this paper, we restrict ourselves to considering only the potential $V$ the same as in Section 4 of \cite{IProbe}.

$\quad$

{\bf\noindent Assumption on $V$.}
Assum that $V(x)$ is a perturbation from the potential that is identically $0$ on a subset of $\Omega$, denoted by $D$, that is, the support of $V$ satisfies
\begin{equation}
%$$
\displaystyle
\text{supp}\,V\subset\overline{D},
\tag {1.2}
%$$
\end{equation}
where $D$ is an open subset  of $\Bbb R^3$  with Lipschitz boundary and satisfies $\overline{D}\subset\Omega$.

$\quad$

\noindent
Furthermore,  we introduce a condition ensuring that $V$ near $\partial D$ exhibits a sharp jump from the background potential, which is identically $0$.

$\quad$

{\bf\noindent Jump condition.}
Let $a\in\partial D$.  We say that $V$ has a {\it positive/negative jump} on $\partial D$ at the point $a$ if there exist a positive constant $C$ and
an open ball $B$ centered at $a$ such that $\pm V(x)\ge C$ for a.e.$x\in B\cap D$.

$\quad$

The Probe Method is a method of probing inside a given domain, by using special input based on a special curve and singular solution of the governing equation
for the background medium.

$\quad$

{\bf\noindent Needle.}  Given $x\in\Omega$ we denote by $N_x$ the set of all piecewise linear curves $\sigma:[0,\,1]\longmapsto\,\overline{\Omega}$ such
that $\sigma$ is one to one as a map, $\sigma(0)\in\partial\Omega$, $\sigma(1)=x$ and $\sigma(t)\in\Omega$ for all $0<t<1$.
Each member $\sigma\in N_x$ is called a {\it needle} with a tip at $x$.

$\quad$

{\bf\noindent Singular Solution.}  
The form of $V$ given by  (1.2) means that the governing equation for the background medium is the Laplace equation.
Let $G(z)$ denote the standard fundamental solution of the Laplace equation, that is
$$\begin{array}{ll}
\displaystyle
G(z)=\frac{1}{4\pi\vert z\vert}, & z\in\mathbb R^3\setminus\{0\}.
\end{array}
$$
This satisfies
$$\displaystyle
\Delta G(z)+\delta(z)=0
$$
in the sense of distribution in the whole space, where $\delta(z)$ denotes the Dirac delta function.
We see that, for each $j=1,2,3$ the function
$$\begin{array}{ll}
\displaystyle
\frac{\partial G}{\partial z_j}(z), & z\in\mathbb R^3\setminus\{0\}, 
\end{array}
$$
satisfies also the Laplace equation in $\mathbb R^3\setminus\{0\}$
and its singularity at $z=0$ is stronger than that of $G(z)$.   
In \cite{IProbe}, we made use of the property: for any finite cone $V$ with the vertex at $0$, 
$$
\displaystyle
\sum_{j=1}^3\int_V\left\vert\frac{\partial G}{\partial z_j}(z)\right\vert^2 dz=\infty.
$$
See Lemma A in the Appendix, which ensures that each term of the integrals on this left-hand side has also the same blowing up property.

$\quad$

Using this singular solution and needle we introduce a sequence of special solutions of the governing equation for the background medium.
The selected sequence of special solutions which are called a {\it needle sequence} for $(x,\sigma)$, is based on Proposition 4 in \cite{IProbe}.  It is an application
of the Runge approximation theorem and stated as below.

$\quad$

{\bf\noindent Needle Sequence.}
Let $j=1,2,3$.
Given $x\in\Omega$ and $\sigma\in N_x$, there exists a sequence of  $\{v_n^j\}$ of harmonic functions in $H^1(\Omega)$ such that,
$\text{supp}\,(v_n^j\vert_{\partial\Omega})\subset\Gamma_0$ and
as $n\rightarrow\infty$
$$\displaystyle
v_n^j\rightarrow\frac{\partial G}{\partial z_j}(\,\cdot\,-x) 
$$
in $H^1_{\text{loc}}(\Omega\setminus\sigma)$, where we write $\sigma([0,\,1])$ as simply $\sigma$ and $\Gamma_0\not=\emptyset$
is an arbitrary fixed open subset of $\partial\Omega$.   We call the sequence $\{v_n^j\}$ the {\it needle sequence} for $(x,\sigma)$.

$\quad$

Using the reformulated Probe Method  \cite{INew}, the relevant part of Section 4 in \cite{IProbe} can be extracted as follows:

\begin{thm}
%\proclaim{\noindent Theorem 1.1.}

\noindent
{\rm (a)}  Let $x\in\Omega\setminus\overline{D}$ and $\sigma\in N_x$ satisfy $\sigma\cap\overline{D}=\emptyset$.  Then, for an arbitrary needle sequence
$\{v_n^j\}$ with $j=1,2,3$ for $(x,\sigma)$ we have
\begin{equation}
%$$\displaystyle
\lim_{n\rightarrow\infty}\sum_{j}<(\Lambda_V-\Lambda_0)v_n^j\vert_{\partial\Omega},v_n^j\vert_{\partial\Omega}>={\mathcal I}(x),
\tag {1.3}
%$$
\end{equation}
where the function ${\mathcal I}(x)$ of independent variable $x\in\Omega\setminus\overline{D}$ is defined by
\begin{equation}
%$$\displaystyle
{\mathcal I}(x)=\int_DV(z)\vert\nabla G(z-x)\vert^2\,dz+\int_{D}V(z)\mbox{\boldmath $w$}_x(z)\cdot\nabla G(z-x)\,dz
\tag {1.4}
%$$
\end{equation}
and the vector valued function $\mbox{\boldmath $w$}_x=\mbox{\boldmath $w$}(z)$ solves
\begin{equation}
%$$
\left\{
\begin{array}{ll}
\displaystyle
-\Delta \mbox{\boldmath $w$}+V(z)\mbox{\boldmath $w$}=-V(z)\nabla G(z-x), & z\in\Omega,
\\
\\
\displaystyle
\mbox{\boldmath $w$}=\mbox{\boldmath $0$}, & z\in\partial\Omega.
\end{array}
\right.
\tag {1.5}
%$$
\end{equation}

\noindent
{\rm (b)}  Let $a$ be an arbitrary point on $\partial D$.  We have
$$\displaystyle
\lim_{x\rightarrow a}{\mathcal I}(x)=
\left\{
\begin{array}{ll}
\displaystyle
\infty & \text{if $V$ has a positive jump on $\partial D$ at point $a$,}\\
\\
\displaystyle
-\infty & \text{if $V$ has a negative jump on $\partial D$ at point $a$.}
\end{array}
\right.
$$

\noindent
{\rm (c)}  For each $\epsilon>0$, it holds that
$$\displaystyle
\sup_{x\in\Omega\setminus\overline{D}, \text{dist}\,(x,\partial D)>\epsilon}\vert{\mathcal I}(x)\vert<\infty.
$$

%\endproclaim
\end{thm}

\noindent
Some remarks on Theorem 1.1 are in order.

$\quad$

{\bf\noindent Remark 1.1.}
Note that, Proposition 7 in Section 4 of \cite{IProbe} can be also restated as follows:
given $a\in\partial D$ one can choose positive constants $\epsilon_0$, $C_1$ and $C_2$ depending on $a$ in such a way that, 
for all $\epsilon\in\,]0,\,\epsilon_0]$, it holds that
$$\begin{array}{ll}
\displaystyle
\vert {\mathcal I}(x)\vert>C_1
\int_{D\cap B_{\epsilon}(a)}\vert\nabla G(z-x)\vert^2 dz-C_2, & \vert x-a\vert<\frac{\epsilon}{2},
\end{array}
$$
where $B_{\epsilon}(a)$ denotes the open ball with radius $\epsilon$ and centered at $a$.
Besides, by the proof of Theorem C in \cite{IProbe} one has the lower estimate as $x\rightarrow a$
$$\displaystyle
\int_{D\cap B_{\epsilon}(a)}\vert\nabla G(z-x)\vert^2 dz\ge\frac{C_3}{\vert x-a\vert},
$$
where $C_3$ is a positive constant.  Therefore we have a lower estimate of the blowing-up rate of $\vert {\mathcal I}(x)\vert$
as $x\rightarrow a$, that is $\vert x-a\vert^{-1}$.

Note that, more precisely, from the proof therein we have: if $V$ has a positive/negative jump on $\partial D$ at point $a$
$$\begin{array}{ll}
\displaystyle
\pm{\mathcal I}(x)>C_1\int_{D\cap B_{\epsilon}(a)}\vert\nabla G(z-x)\vert^2 dz-C_2, & \vert x-a\vert<\frac{\epsilon}{2}.
\end{array}
$$
Besides, Theorem E  in Section 5 of \cite{IProbe} clarified the behaviour of ${\mathcal I}(x)$ as $x\rightarrow a$:
$$\begin{array}{ll}
\displaystyle
{\mathcal I}(x)\sim \alpha\int_{D}\vert\nabla G(z-x)\vert^2\,dz, & x\rightarrow a,
\end{array}
$$
where $\alpha=\lim_{D\ni x\rightarrow a}V(x)$ provided its existence and $\alpha\not=0$.
Note that, not mentioned there, from this one can deduce also its boundary local version
$$\begin{array}{ll}
\displaystyle
{\mathcal I}(x)\sim \alpha\int_{D\cap B_{\epsilon}(a)}\vert\nabla G(z-x)\vert^2\,dz, & x\rightarrow a,
\end{array}
$$
since
$$\begin{array}{ll}
\displaystyle
\int_D\vert\nabla G(z-x)\vert^2\,dz\sim\int_{D\cap B_{\epsilon}a)}\vert\nabla G(z-x)\vert^2\,dz, & x\rightarrow a.
\end{array}
$$
The formula gives us an explicit analytical formula for the calculation of the jump of the potential across $\partial D$
provided $D$ is {\it known}.

$\quad$

{\bf\noindent Remark 1.2.}
Given $x\in\Omega\setminus\overline{D}$ let $\sigma\in N_x$ satisfy $\sigma\cap\overline{D}=\emptyset$.
Let $\{v_n\}$ be the original needle sequence for $(x,\sigma)$ in \cite{INew}, that is, a sequence of solutions of the Laplace equation in $H^1(\Omega)$ such that
$$\displaystyle
v_n\rightarrow G(\,\cdot\,-x) 
$$
in $H^1_{\text{loc}}(\Omega\setminus\sigma)$ (see also Proposition 1 in \cite{IProbe} for its origin).
Then, a similar argument for the derivation of (1.3) yields also
$$\displaystyle
\lim_{n\rightarrow\infty}<(\Lambda_V-\Lambda_0)v_n\vert_{\partial\Omega},v_n\vert_{\partial\Omega}>=I(x),
$$
where
$$\displaystyle
I(x)=\int_D V(z)\vert G(z-x)\vert^2\,dz+V(z)w_x(z)G(z-x)\,dz
$$
and $w_x=w$ solves
$$\left\{
\begin{array}{ll}
\displaystyle
-\Delta w+V(z)w=-V(z)G(z-x), & z\in\Omega,
\\
\\
\displaystyle
w=0, & z\in\partial\Omega.
\end{array}
\right.
$$
Here we have, for a positive constant $C$ independent of $x\in\Omega$
$$\displaystyle
\Vert w_x\Vert_{L^2(\Omega)}\le C\Vert V\Vert_{L^{\infty}(\Omega)}\Vert G(\,\cdot\,-x)\Vert_{L^2(\Omega)}.
$$
Besides it holds that, for $\alpha\in\,]0,\,3[$
$$\displaystyle
\sup_{x\in\Bbb R^3}\,\int_{\Omega}\frac{dz}{\vert x-z\vert^{\alpha}}<\infty.
$$
Thus one gets
$$\displaystyle
\sup_{x\in\Omega\setminus\overline{D}}\,\vert I(x)\vert<\infty.
$$
This means that the $I(x)$ does not play the role of the indicator function in the Probe Method for the equation (1.1).
Clearly, this is because of the weak singularity of $G(z)$ itself.
This is the main and essential reason why we used the needle sequences generated from the derivatives of $G(z)$.
This is the decisive difference from the penetrable obstacle case appearing as a jump of second-order coefficient of the governing equation, see Proposition 1 of \cite{IProbe}.

\subsection {IPS and Side A of Singular Sources Method}

Now we describe the main object of this paper.  As can be seen from (a) to (c) of Theorem 1.1, we have the {\it indicator function} for the Probe Method, that is, ${\mathcal I}(x)$
and the {\it reflected solution} should be the $\mbox{\boldmath $w$}_x$ that is the solution of (1.5).
Considering the spirit of IPS, a candidate of the indicator function for the Singular Sources Method should be
$$\begin{array}{l}
\Omega\setminus\overline{D}\ni x\longmapsto \mbox{\boldmath $w$}_x(x),
\end{array}
$$
where $\mbox{\boldmath $w$}_x(x)=\mbox{\boldmath $w$}_x(z)\vert_{z=x}$.
However, note that the previous IPS developed in \cite{IPS}, \cite{IPS2} and \cite{IPS3} is based on the needle sequences and reflected solutions which are
generated by the {\it fundamental solutions} themself not their derivatives of the governing equations for the background media.
Thus we raise the problem as follows.

$\quad$

{\bf\noindent Problem 1.}  Find the IPS function between ${\mathcal I}(x)$ and $\mbox{\boldmath $w$}_x(x)$.  In other words, develop the IPS by using
${\mathcal I}(x)$ and $\mbox{\boldmath $w$}_x$ and clarify what is the Singular Sources Method in the framework of IPS.

$\quad$

It should be emphasized that the $\mbox{\boldmath $w$}_x(x)$  is the vector valued function in contrast to the scalar valued function ${\mathcal I}(x)$.
Thus the behaviour of $\mbox{\boldmath $w$}_x(x)$ as $x$ approaches a point on $\partial D$ would not be simple.  This is also a new situation.

In what follows we denote by $\nabla \mbox{\boldmath $u$}$ the Jacobian matrix of the vector valued function $\mbox{\boldmath $u$}$.
Thus the scalar field $\text{Trace}\,\nabla\mbox{\boldmath $u$}$ is the divergence of $\mbox{\boldmath $u$}$ and denoted by $\nabla\cdot\mbox{\boldmath $u$}$.
For $3\times3$ matrixes $A=(a_{ij})$ and $B=(b_{ij})$ we write $A\cdot B=\sum_{ij}a_{ij}b_{ij}$ and $\vert A\vert=\sqrt{A\cdot A}$.

Our main discovery is as follows.

\begin{thm}
%\proclaim{\noindent Theorem 1.2.}
We have
\begin{equation}
%$$
\begin{array}{ll}
\displaystyle
(\nabla\cdot\mbox{\boldmath $w$}_x)(x)+(\nabla\cdot\mbox{\boldmath $w$}_x^1)(x)
={\mathcal I}(x)+{\mathcal I}^1(x), & x\in\Omega\setminus\overline{D},
\end{array}
\tag {1.6}
%$$
\end{equation}
where the vector valued function $\mbox{\boldmath $w$}_x^1=\mbox{\boldmath $w$}(z)$ solves 
\begin{equation}
%$$
\left\{
\begin{array}{ll}
\displaystyle
-\Delta \mbox{\boldmath $w$}+V(z)\mbox{\boldmath $w$}=\mbox{\boldmath $0$}, & z\in\Omega,
\\
\\
\displaystyle
\mbox{\boldmath $w$}=\nabla G(z-x), & z\in\partial\Omega
\end{array}
\right.
\tag {1.7}
%$$
\end{equation}
and function ${\mathcal I}^1(x)$ of independent variable $x\in\Omega\setminus\overline{D}$  is defined by
\begin{equation}
%$$
\displaystyle
{\mathcal I}^1(x)=-\int_{\Omega}\vert\nabla\mbox{\boldmath $w$}_x^1\vert^2\,dz
-\int_{D}V(z)\vert\mbox{\boldmath $w$}_x^1\vert^2\,dz
+\int_{\partial\Omega}\nabla G(z-x)\cdot\frac{\partial}{\partial\nu}\nabla G(z-x)\,dS(z).
\tag {1.8}
%$$
\end{equation}

%\endproclaim

\end{thm}

$\quad$

{\bf\noindent Remark 1.3.}
As a direct consequence of (2.13) and (2.24) in Subsection 2.2, we obtain more impressive representation:
\begin{equation}
%$$
\displaystyle
{\mathcal I}(x)=-\int_{\Omega}\vert\nabla\mbox{\boldmath $w$}_x\vert^2\,dz
-\int_{D}V(z)\vert\mbox{\boldmath $w$}_x\vert^2\,dz+\int_{D}V(z)\vert\nabla G(z-x)\vert^2\,dz;
\tag {1.9}
%$$
\end{equation}
\begin{equation}
%$$
\displaystyle
{\mathcal I}^1(x)=-\int_{\Omega}\vert\nabla\mbox{\boldmath $w$}_x^1\vert^2\,dz
-\int_{D}V(z)\vert\mbox{\boldmath $w$}_x^1\vert^2\,dz
-\int_{\Bbb R^3\setminus\overline{\Omega}}\vert\nabla^2 G(z-x)\vert^2\,dz.
\tag {1.10}
%$$
\end{equation}
Here $\nabla^2 G$ denotes the Jacobian matrix of $\nabla G$.
Expressions (1.9) and (1.10) also strongly support that the right-hand side term of equation (1.6) should be called the inner decomposition of IPS function; 
see \cite{IPS} for its origin.

$\quad$

From a pont of the IPS view, this suggests, the candidate of the indicator function of the Singular Sources Method via the Probe Method,
as stated in Theorem 1.1,  should be
the scalar valued function defined by
\begin{equation}
%$$
\displaystyle
\Omega\setminus\overline{D}\ni x\longmapsto (\nabla\cdot\mbox{\boldmath $w$}_x)(x).
\tag {1.11}
%$$
\end{equation}
In order to to justify that this is the indicator function for the Singular Sources Method in IPS, we
have to give a calculation  method of the value of the function (1.11) at an arbitrary point in $\Omega\setminus\overline{D}$.
For this purpose, we introduce a special sequence of functions which plays the key role.

$\quad$

{\bf\noindent Definition 1.1.}  Given $x\in\Omega$ and $\sigma\in N_x$, let $\{v_n^j\}$ be a needle sequence for $(x,\sigma)$.
Define
\begin{equation}
%$$
\begin{array}{ll}
\displaystyle
G_n^j(z,x)=\frac{\partial G}{\partial z_j}(z-x)-v_n^j(z), 
&
\displaystyle
z\in\Omega\setminus\{x\}.
\end{array}
\tag {1.12}
%$$
\end{equation}
The sequence satisfies as $n\rightarrow\infty$, $G_n^j(\,\cdot\,,x)\rightarrow 0$ in $H^1(B)$ for any open ball $B$ with $\overline{B}\subset\Omega\setminus\sigma$.

$\quad$

The next theorem shows that the function defined by (1.11) plays the role of the indicator function in the Singular Sources Method in IPS.

\begin{thm}
%\proclaim{\noindent Theorem 1.3.}

\noindent
{\rm (a)}  Let $x\in\Omega\setminus\overline{D}$ and $\sigma\in N_x$ satisfy $\sigma\cap\overline{D}=\emptyset$.  Then, for an arbitrary needle sequence
$\{v_n^j\}$ with $j=1,2,3$ for $(x,\sigma)$, we have
\begin{equation}
%$$
\displaystyle
-\lim_{n\rightarrow\infty}\sum_{j}<(\Lambda_V-\Lambda_0)v_n^j\vert_{\partial\Omega},G_n^j(\,\cdot\,,x)\vert_{\partial\Omega}>=(\nabla\cdot\mbox{\boldmath $w$}_x)(x),
\tag {1.13}
%$$
\end{equation}

\noindent
{\rm(b)}  For each $\epsilon>0$ we have
\begin{equation}
%$$
\displaystyle
\sup_{x\in\Omega\setminus\overline{D},\,\text{dist}\,(x,\partial\Omega)>\epsilon}
\vert \nabla\cdot\mbox{\boldmath $w$}_x^1(x)\vert<\infty,
\tag {1.14}
%$$
\end{equation}
\begin{equation}
%$$
\displaystyle
\sup_{x\in\Omega\setminus\overline{D},\,\text{dist}\,(x,\partial\Omega)>\epsilon}
\vert{\mathcal I}^1(x)\vert<\infty
\tag {1.15}
%$$
\end{equation}
and
\begin{equation}
%$$
\displaystyle
\sup_{x\in\Omega\setminus\overline{D},\,\text{dist}\,(x,\partial\Omega)>\epsilon}
\vert\nabla\cdot\mbox{\boldmath $w$}_x(x)-{\mathcal I}(x)\vert<\infty.
\tag {1.16}
%$$
\end{equation}

\noindent
{\rm (c)}  Let $a$ be an arbitrary point on $\partial D$.  We have
$$\displaystyle
\lim_{x\rightarrow a}(\nabla\cdot\mbox{\boldmath $w$}_x)(x)=
\left\{
\begin{array}{ll}
\displaystyle
\infty & \text{if $V$ has a positive jump on $\partial D$ at point $a$,}\\
\\
\displaystyle
-\infty & \text{if $V$ has a negative jump on $\partial D$ at point $a$.}
\end{array}
\right.
$$

\noindent
{\rm (d)}  For each $\epsilon_i>0$, $i=1,2$ it holds that
$$\displaystyle
\sup_{x\in\Omega\setminus\overline{D},\,\text{dist}\,(x,\partial D)>\epsilon_1,\,\text{dist}\,(x,\partial\Omega)>\epsilon_2}\,\vert(\nabla\cdot\mbox{\boldmath $w$}_x)(x)\vert<\infty.
$$

%\endproclaim

\end{thm}

Some remarks are in order.

$\quad$

\noindent
$\bullet$ No regularity assumption is assumed except for $V\in L^{\infty}(\Omega)$.

\noindent
$\bullet$  The $V$ can take either positive or negative values in $D$ except for a neighbourhod of $\partial D$ relative to $\overline{D}$.

\noindent
$\bullet$  By using (1.16), all the results on ${\mathcal I}(x)$ mentioned as Remark 1.1 also can be transplanted.

\noindent
$\bullet$  Theorem 1.3 should be called the Side A of Singular Sources Method via the Probe Method in the framework of IPS
or simply, Singular Sources Method in IPS.
In particular, the statement (c) of Theorem 1.3, which is a consequence of {\rm(b)} of Theorem 1.1 and (1.6)
is quite {\it intuitive}.  This gives a characterization of the surface of {\it penetrable obstacle} as the blowing up set
of the restriction of  the {\it divergence} of the vector field $\mbox{\boldmath $w$}_x(y)$ onto $y=x$.

\noindent
$\bullet$  It seems to be impossible to delete in {\rm(d)} the restriction $\text{dist}\,(x,\partial\Omega)>\epsilon_2$.  From the well posedness of 
Dirchet problem (1.5) one can obtain the boundedness of $\Vert\mbox{\boldmath $w$}_x\Vert_{H^2(\Omega)}$ with $\text{dist}(x,\partial D)>\epsilon_1$.
Thus by the Sobolev embedding, one has the boundedness of $\Vert\mbox{\boldmath $w$}_x\Vert_{L^{\infty}(\Omega)}$ with $\text{dist}\,(x,\partial D)>\epsilon_1$.
However, potential $V$ is just essentially bounded, we may not do such an argument for  $\nabla\cdot\mbox{\boldmath $w$}_x$ directly.
This is a new situation different from the previous results for impenetrable obstacles, see
\cite{IPS}, \cite{IPS2} and \cite{IPS3}.

$\quad$

Theorems 1.1 to 1.3 suggest us that we should define the IPS function for the Schr\"odinger equation $-\Delta u+V(z)u=0$ as follows.

$\quad$

{\bf\noindent Definition 1.2.}
Given $x\in\Omega\setminus\overline{D}$ let $\mbox{\boldmath $W$}_x=\mbox{\boldmath $W$}$ solve
$$\left\{
\begin{array}{ll}
\displaystyle
-\Delta \mbox{\boldmath $W$}+V(z)\mbox{\boldmath $W$}=-V(z)\nabla G(z-x), & z\in\Omega,
\\
\\
\displaystyle
\mbox{\boldmath $W$}=\nabla G(z-x), & z\in\partial\Omega.
\end{array}
\right.
$$
We call the function
$$\displaystyle
\Omega\setminus\overline{D}\ni x\longmapsto (\nabla\cdot\mbox{\boldmath $W$}_x)(y)\vert_{y=x},
$$
the {\it IPS function} for the Schr\"odinger equation $-\Delta u+V(x)u=0$.

$\quad$

\noindent
Since we have the expression
$$\begin{array}{ll}
\displaystyle
\mbox{\boldmath $W$}_x(z)=\mbox{\boldmath $w$}_x(z)+\mbox{\boldmath $w$}_x^1(z), & z\in \Omega,
\end{array}
$$
the outer (point-wise) and inner (energy integral) decompositions of IPS function, which is the symbol of
the integration of the Probe and Singular Sources Methods is the formula
$$\displaystyle
\nabla\cdot\mbox{\boldmath $w$}_x(x)+\nabla\cdot\mbox{\boldmath $w$}^1_x(x)=\nabla\cdot\mbox{\boldmath $W$}_x(x)=
{\mathcal I}(x)+{\mathcal I}^1(x).
$$
Besides, by (1.14), (1.15) one has, as $x\rightarrow a\in\partial D$
we have
$$\displaystyle
\nabla\cdot\mbox{\boldmath $w$}_x(x)\sim\nabla\cdot\mbox{\boldmath $W$}_x(x)\sim {\mathcal I}(x)
$$
modulo $O(1)$.  Thus it follows from {\rm (b)} of Theorem 1.1 that, for each $a\in\partial D$ we have
$$\displaystyle
\lim_{x\rightarrow a}\,\nabla\cdot\mbox{\boldmath $W$}_x(x)=
\left\{\begin{array}{ll}
\displaystyle
\infty & \text{if $V$ has a positive jump on $\partial D$ at $a$,}
\\
\\
\displaystyle
-\infty & \text{if $V$ has a negative jump on $\partial D$ at $a$.}
\end{array}
\right.
$$
Furthermore it holds that, for each point $b\in\partial\Omega$
$$\displaystyle
\lim_{x\rightarrow b}\nabla\cdot\mbox{\boldmath $W$}_x(x)=-\infty.
$$
This needs some technical argument and for the proof see Appendix.

Here we point out an expression of the IPS function $\nabla\cdot\mbox{\boldmath $W$}_x(x)$ in terms of energy integral.
That is,
$$\displaystyle
\nabla\cdot\mbox{\boldmath $W$}_x(x)=-\int_{\Omega}\vert\nabla\mbox{\boldmath $W$}_x\vert^2+V(z)\vert\mbox{\boldmath $W$}_x\vert^2\,dz
+\int_{D}V(z)\vert\nabla G(z-x)\vert^2\,dz-\int_{\Bbb R^3\setminus\overline{\Omega}}\vert\nabla^2G(z-x)\vert^2\,dz.
$$
This expression corresponds to that of the IPS function in an impenetrable obstacle case, see Remark 1.7 of \cite{IPS}.
This is a direct corollary of the following fact which is proved in Subsection 2.2.

$\quad$

{\bf\noindent Proposition 1.1.}  {\it It holds that, for all $(x,y)\in(\Omega\setminus\overline{D})^2$}
$$\begin{array}{l}
\displaystyle
\,\,\,\,\,\,
\frac{\nabla\cdot\mbox{\boldmath $W$}_x(y)+\nabla\cdot\mbox{\boldmath $W$}_y(x)}{2}+\int_{\Omega}\nabla\mbox{\boldmath $W$}_x\cdot\nabla\mbox{\boldmath $W$}_y+V(z)\mbox{\boldmath $W$}_x\cdot\mbox{\boldmath $W$}_y\,dz
\\
\\
\displaystyle
=
\int_{\Omega} V(z)\nabla G(z-x)\cdot\nabla G(z-y)\,dz
-\int_{\Bbb R^3\setminus\overline{\Omega}}\nabla^2 G(z-x)\cdot\nabla^2 G(z-y)\,dz.
\end{array}
$$

\noindent
This can be considered as an integro-differential equation and the family $(\mbox{\boldmath $W$}_{\zeta})_{\zeta\in\Omega\setminus\overline{D}}$ of $H^2(\Omega)$
is an exact solution.
Note that, in \cite{IPS},  we have already seen a similar integro-differential equation for an inverse obstacle problem governed by the Laplace equation
with homogeneous Neumann boundary condition on the surface of unknown obstacle, see Remark 3.1 of \cite{IPS}.

Therefore from these evidence, it is reasonable to refer to the function $\nabla\cdot\mbox{\boldmath $W$}_x(x)$ the IPS function for the equation $-\Delta u+V(x)u=0$.

\subsection{Completely Integrated IPS}

In this subsection, we present a completely integrated method that combines both the Probe Method and Singular Sources Method.

Given $x\in\Omega$ and $\sigma\in N_x$ let $\{v_n^j\}$, $j=1,2,3$ be an arbitrary needle sequence for $(x,\sigma)$.
First we introduce a modification of the needle sequence $\{v_n^j\}$ for $(x,\sigma)$.

\noindent
Let the vector valued function $\mbox{\boldmath $H$}(\,\cdot\,,x)$ be the solution of
$$\displaystyle
\left\{
\begin{array}{ll}
\displaystyle
\Delta \mbox{\boldmath $H$}=\mbox{\boldmath $0$}, & z\in\Omega,
\\
\\
\displaystyle
\mbox{\boldmath $H$}=-\nabla G(z-x), & z\in\partial\Omega.
\end{array}
\right.
$$
Note that the minus sign of the Dirichlet data is important.
We denote by $H^j(\,\cdot\,,x)$ the $j$-component of $\mbox{\boldmath $H$}(\,\cdot\,,x)$.
Clearly the sequence $\{v_n^j+H^j(\,\cdot\,,x)\}$ is also a needle sequence for $(x,\sigma)$ in the sense that
\begin{equation}
%$$
\displaystyle
v_n^j(\,\cdot\,,x)+H^j(\,\cdot\,,x)\rightarrow \frac{\partial G}{\partial z_j}(\,\cdot\,-x)+H^j(\,\cdot\,,x)\equiv G^j(\,\cdot\,,x)
\tag {1.17}
%$$
\end{equation}
in $H^1_{\text{loc}}(\Omega\setminus\sigma))$.  Note that the singularity of $G^j(z,x)$ at $z=x$ coincides with that of $\frac{\partial G}{\partial z_j}(z-x)$.

Then for $j=1,2,3$ the trace of $G_n^j(\,\cdot\,,x)$ defined by (1.12) onto $\partial\Omega$ has  a tricky expression
\begin{equation}
%$$
\begin{array}{l}
\displaystyle
G_n^j(\,\cdot\,,x)\vert_{\partial\Omega}=-(v_n^j+H^j(\,\cdot\,,x))\vert_{\partial\Omega},
\end{array}
\tag{1.18}
%$$
\end{equation}
This means that the harmonic function $-(v_n^j+H^j(\,\cdot\,,x))$ is a lifting of function $G_n^j(\,\cdot\,,x)\vert_{\partial\Omega}$ into the whole domain $\Omega$.

The following theorem provides an exact integration of the Probe Method and Singular Sources Method.

\begin{thm}
%\proclaim{\noindent Theorem 1.4.}

\noindent
{\rm (a)}  Let $x\in\Omega\setminus\overline{D}$ and $\sigma\in N_x$ satisfy $\sigma\cap\overline{D}=\emptyset$.  Then, for an arbitrary needle sequence
$\{v_n^j\}$ with $j=1,2,3$ for $(x,\sigma)$ we have
\begin{equation}
%$$
\displaystyle
\lim_{n\rightarrow\infty}\sum_{j}<(\Lambda_V-\Lambda_0)G_n^j(\,\cdot\,,x)\vert_{\partial\Omega},
G_n^j(\,\cdot\,,x))\vert_{\partial\Omega}>={\mathcal I}^*(x),
\tag {1.19}
%$$
\end{equation}
where the function ${\mathcal I}^*(x)$ of independent variable $x\in\Omega\setminus\overline{D}$ is defined by
\begin{equation}
%$$
\displaystyle
{\mathcal I}^*(x)=
\int_DV(z)\vert\nabla G(z-x)+\mbox{\boldmath $H$}(z,x)\vert^2\,dz+\int_{D}V(z)\mbox{\boldmath $w$}^*_x(z)\cdot
\left(\nabla G(z-x)+\mbox{\boldmath $H$}(z,x)\right)\,dz
\tag {1.20}
%$$
\end{equation}
and the vector valued function $\mbox{\boldmath $w$}^*_x=\mbox{\boldmath $w$}(z)$ solves
\begin{equation}
%$$
\left\{
\begin{array}{ll}
\displaystyle
-\Delta \mbox{\boldmath $w$}+V(z)\mbox{\boldmath $w$}=-V(z)(\nabla G(z-x)+\mbox{\boldmath $H$}(z,x)), & z\in\Omega,
\\
\\
\displaystyle
\mbox{\boldmath $w$}=\mbox{\boldmath $0$}, & z\in\partial\Omega.
\end{array}
\right.
\tag {1.21}
%$$
\end{equation}

\noindent
{\rm (b)}  It holds that
\begin{equation}
%$$
\displaystyle
{\mathcal I}^*(x)={\mathcal I}(x)+2({\mathcal I}^1(x)-\nabla\cdot\mbox{\boldmath $w$}^1_x(x))+
<(\Lambda_V-\Lambda_0)\nabla G(\,\cdot\,-x)\vert_{\partial\Omega}, \nabla G(\,\cdot\,-x)\vert_{\partial\Omega}>.
\tag {1.22}
%$$
\end{equation}

\noindent
{\rm (c)}  Let $a$ be an arbitrary point on $\partial D$.  We have
$$\displaystyle
\lim_{x\rightarrow a}{\mathcal I}^*(x)=
\left\{
\begin{array}{ll}
\displaystyle
\infty & \text{if $V$ has a positive jump on $\partial D$ at point $a$,}\\
\\
\displaystyle
-\infty & \text{if $V$ has a negative jump on $\partial D$ at point $a$.}
\end{array}
\right.
$$

\noindent
{\rm (d)}  For each $\epsilon_i>0$, $i=1,2$ it holds that
$$\displaystyle
\sup_{x\in\Omega\setminus\overline{D}, \text{dist}\,(x,\partial D)>\epsilon_1, \text{dist}\,(x,\partial\Omega)>\epsilon_2}\vert{\mathcal I}^*(x)\vert<\infty.
$$

\noindent
{\rm (e)}  It holds that
\begin{equation}
%$$
\begin{array}{ll}
\displaystyle
(\nabla\cdot\mbox{\boldmath $w$}^*_x)(x)
={\mathcal I}^*(x), & x\in\Omega\setminus\overline{D}.
\end{array}
\tag {1.23}
%$$
\end{equation}

\noindent
{\rm (f)}  It holds that
\begin{equation}
%$$
\displaystyle
\nabla\cdot\mbox{\boldmath $w$}^*_x(x)=\nabla\cdot\mbox{\boldmath $w$}_x(x)
+({\mathcal I}^1(x)-\nabla\cdot\mbox{\boldmath $w$}^1_x(x))
+<(\Lambda_V-\Lambda_0)\nabla G(\,\cdot\,,x)\vert_{\partial\Omega},\nabla G(\,\cdot\,-x)\vert_{\partial\Omega}>.
\tag {1.24}
%$$
\end{equation}

%\endproclaim
\end{thm}

Some remarks are in order.

$\quad$

\noindent
$\bullet$  
The formula  (1.23) justifies the name ``completely integrated IPS'' for the method based on Theorem 1.4.

\noindent
$\bullet$  The function ${\mathcal I}^*$ defined by (1.20) and the governing equation (1.21) of vector valued function $\mbox{\boldmath $w$}^*_x$ are described
by using the vector valued harmonic function $\mbox{\boldmath $H$}(\,\cdot\,,x)$.  However, in the formulae (1.22) and (1.24) this function never appears.  Those formulae represent a relationship between the completely integrated version of IPS and original one.
Those correspond to formulae (4.8) of Theorem 7 and (4.13) of Remark 4.2 in \cite{IPS}.

$\quad$

\subsection{Comparison to Singular Sources Method of Potthast}

In \cite{Pbook} Potthast considered an inverse medium scattering problem of acoustic wave and introduced
the Singular Sources Method using {\it multipoles of higher order}, which are singular solutions of the Helmholtz equation governing the background acoustic medium.
First let us review his method in three dimensions by taking up
the most typical case, since it is closely related to Theorem 1.3.

He employed the point source $\Phi_{1,q}(\,\cdot\,,x)$ of multipole of order one (p.66 in \cite{Pbook}) at $x\in\Bbb R^3$, that is
$$\begin{array}{ll}
\displaystyle
\Phi_{1,q}(z,x)=\frac{1}{4\pi}h^{(1)}_1(k\vert z-x\vert)\frac{(z-x)\cdot q}{\vert z-x\vert}, & z\not=x.
\end{array}
$$
where $k>0$, $q\in S^2$ denotes the polarization and $h_1^{(1)}$ the spherical Hankel function of the first kind.
Since 
$$\displaystyle
h^{(1)}_1(\zeta)=-e^{ik\zeta}\left(\frac{1}{\zeta^2}+\frac{i}{\zeta}\right),
$$
we have the expression
$$\displaystyle
\Phi_{1,q}(z,x)=-\nabla\Phi(z-x)\cdot q,
$$
where
$$\begin{array}{ll}
\displaystyle
\Phi(z)=\frac{e^{ik\vert z\vert}}{4\pi\vert z\vert}, & z\in\Bbb R^3\setminus\{x\}.
\end{array}
$$
Thus $\Phi_{1,q}(\,\cdot\,,x)$ is a singular solution of the Helmholtz equation
$$\begin{array}{ll}
\displaystyle
\Delta v+k^2 v=0, & z\in\Bbb R^3\setminus\{x\}.
\end{array}
$$
He considered the scattered field denoted by $\Phi^s_{1,q}(\,\cdot\,,x)=w$ for incident multipole $\Phi_{1,q}(\,\cdot\,,x)$, which is the unique solution of
inhomogeneous Helmholtz equation
\begin{equation}
%$$
\begin{array}{ll}
\displaystyle
\Delta w+k^2 n(z)w=-k^2(n(z)-1)\Phi_{1,q}(z,x), & z\in\Bbb R^3
\end{array}
\tag {1.25}
%$$
\end{equation}
with the Sommerfeld radiation condition.

The assumption on $n(z)$ (Definition 2.2.5 on p.62 in \cite{Pbook}) is as follows:

\noindent
The support condition
$$\displaystyle
\text{supp}\,(n-1)\subset\overline{D}.
$$

\noindent
The regularity condition: $\partial D\in C^2$ and for some $\alpha\in\,[0,\,1]$
$$\displaystyle
n\in C^{0,\alpha}(\overline{D}).
$$

\noindent
The jump condition:
there exist positive constants $C_1$ and $C_2$ such that, for all $z$ in a neighbourhood of $\partial D$ relative to $\overline{D}$
$$\displaystyle
C_1\le \vert n(z)-1\vert\le C_2.
$$
The Singular Sources Method for inverse medium scattering problem makes use of the field
$$\begin{array}{ll}
\displaystyle
\Phi^s_{1,q}(z,x), & (z,x)\in(\Bbb R^3\setminus\overline{D})^2.
\end{array}
$$
Note that this field depends also on polarization $q$.

Roughly speaking, his method consists of two steps:

$\quad$

\noindent
{\bf\noindent Step 1.}
 Calculation formula of the scattered field for incident multipole of order one from the {\it far field operator}, 
which is the observation data (Theorem 3.1.6 on p.131 of \cite{Pbook}).

$\quad$

\noindent
{\bf\noindent Step 2.}
Blowing up of the value of scattered field at the source point when source point approaches the surface of obstacle (Theorem 2.1.12 on p.73 of \cite{Pbook}).
More precisely he proved that: given $a\in\partial D$, as $h\rightarrow 0$
$$\displaystyle
\vert \Phi^s_{1,q}(a+\nu(a)h,a+\nu(a)h)\vert_{q=-\nu(a)}\vert\equiv \vert\Phi^s_{1,-\nu(a)}(a+\nu(a)h,a+\nu(a)h)\vert\rightarrow\infty,
$$
where $\nu(a)$ denotes the unit outward normal to the surface at point $a$.  
The blowing-up rate given therein is $\vert\log h\vert$ as $h\rightarrow 0$.

$\quad$

Some comments are in order.
Step 1 makes use of his point source method in \cite{P0} (1998) which is an analytical continuation procedure of scattered field from its far field pattern.
The result in Step 2 says that 
when $x$ approaches the surface along the line $x=a+\nu(a)h$, with $0<h<<1$, by choosing polarization $q=-\nu(a)$, one gets the blowing-up
of the field $\Phi^s_{1,q}(x,x)$.  This means one needs the knowledge
about  the normal $\nu(a)$ in advance.
There is no description about the behaviour of $\vert\Phi^s_{1,q}(x,x)\vert$  not only as $x$ simply approaches the point $a$ 
but also for an arbitrary {\it fixed} polarization $q$.
This means that the treatment of polarization is not complete from a mathematical point of view
(see also the sentence starting with ``From a practical point of view...'' on pages 204-205 of \cite{Pbook} for his comment on this point).

Hereafter let us explain the difference of our approach from his one.
By the superposition, we see that his field $\Phi^s_{1,q}(\,\cdot\,,x)$ has the expression
\begin{equation}
%$$
\displaystyle
\Phi^s_{1,q}(z,x)=\sum_{j}q_j\Phi^s_{1,\mbox{\boldmath $e$}_j}(z,x),
\tag {1.26}
%$$
\end{equation}
where $q_j$ denotes the $j$-th component of polarization vector $q$ and $\mbox{\boldmath $e$}_j$ the $j$-th standard basis vector of $\Bbb R^3$.
Here, let $(\mbox{\boldmath $w$}_x)_j$ denote the $j$-th component of $\mbox{\boldmath $w$}_x$ which is the solution of (1.5).
Comparing (1.25) to the $j$-th component of  (1.5), we see that field $\Phi^s_{1,\mbox{\boldmath $e$}_j}(\,\cdot\,,x)$ corresponds to $(\mbox{\boldmath $w$}_x)_j$.
By this correspondence, if one considers an analogy in our problem of his field $\Phi^s_{1,q}(z,x)$ expressed by (1.26),  it would be natural to employ the scalar field defined by
$$\displaystyle
w(z)=\sum_{j}q_j(\mbox{\boldmath $w$}_{x})_j(z),
$$
which solves
$$\left\{
\begin{array}{ll}
\displaystyle
-\Delta w+V(z)w=-V(z)\nabla G(z-x)\cdot q, & z\in\Omega,
\\
\\
\displaystyle
w=0, & z\in\partial\Omega.
\end{array}
\right.
$$
However, against this, Theorem 1.2,  which is the core of IPS, suggests us to employ the field
\begin{equation}
%$$
\displaystyle
\nabla\cdot\mbox{\boldmath $w$}_x(z)=\sum_{j}\frac{\partial (\mbox{\boldmath $w$}_x)_j}{\partial z_j}(z).
\tag {1.27}
%$$
\end{equation}
As Theorem 1.3 demonstrates this is right.

$\quad$

\noindent
$\bullet$  The choice (1.27) avoids a technical difficulty caused by the presence of polarization.
One obtains the blowing-up property of $\nabla\cdot\mbox{\boldmath $w$}_x(x)$ as $x$ simply approaches point $a\in\partial D$
without considering polarization.  Besides, by using (1.16) and Remark 1.1, one obtains the lower estimate 
of the bowing-up rate as $x\rightarrow a$
$$\displaystyle
\pm\nabla\cdot\mbox{\boldmath $w$}_x(x)\ge\frac{C_4}{\vert x-a\vert},
$$
where $C_4$ is a positive constant.
This is better than the blowing-up rate described in Step 2.  Note that, in our case $\partial D$ is just Lipschitz \cite{Gr}.

$\quad$

\noindent
These show the effectivity of IPS, and {\it decisive difference} in the concept of the Singular Sources Method within IPS framework
compared to Potthast's original method.

$$\quad$$

This paper is organized as follows.  Theorem 1.2 is proved in Subsection 2.1 as a special case of Proposition 2.1. 
 It is derived by using Green's representation formula  of solutions of (1.5) and (1.7), that is, $\mbox{\boldmath $w$}_x$ and $\mbox{\boldmath $w$}_x^1$.
This proposition introduces a decomposition formula of $\nabla\cdot\mbox{\boldmath $W$}_x(y)$ for all $(x,y)\in(\Omega\setminus\overline{D})^2$ in terms of
the functions ${\mathcal I}(x,y)$ and ${\mathcal I}^1(x,y)$ defined by (2.2) and (2.3), which we call the liftings of ${\mathcal I}(x)$ and ${\mathcal I}^1(x)$, respectively.
Subsection 2.2 lists and proves their properties, especially the energy integral representation, symmetry, and harmonicity when one of the two variables $x$ and $y$ is fixed.
Theorems 1.3 and  1.4 are proved in Subsections 2.3 and 2.4, respectively.  
In Subsection 2.5, as an advantage of introducing the lifting of ${\mathcal I}(x)$,
we describe Theorem 2.1 which states that the local proving in an arbitrary (small) neighbourhood of a given point on $\partial\Omega$ relative to $\overline{\Omega}$
uniquely determines the indicator function ${\mathcal I}(x)$ in the domain $\Omega\setminus\overline{D}$ globally.
Section 3 is devoted to Side B of IPS, which states the blowing up property of
 three types of sequences calculated from Dirichlet-to-Neumann map.   They are called the indicator sequences and in Theorem 3.1, under stronger assumption on
 the jump of  $V$,  their blowing up are proved.  The key for the proof is Proposition 3.2, which states, roughly speaking, on an open set
 having intersection with a needle, the corresponding needle sequence blows up and  its $L^1$-norm is weaker than $L^2$-norm.
 This last property is a new fact.
 As a consequence of Theorem 3.1 a characterization of obstacle is derived as Corollaries 3.1 to 3.3.
 The final section describes some remarks on problems for future consideration.  
 In the Appendix, first we provide the proof of Lemma A for the reader's convenience.
 This lemma is needed for Proposition 3.1 of Section 3, which concernes the blowing up property of needle sequences.
 The second is devoted to a proof of the blowing-up of the IPS function on $\partial\Omega$.

\section{Side A of IPS}

\subsection{Proof of Theorem 1.2}

$\quad$

Theorem 1.2 can be considered as a restricted version of the following statement.

$\quad$

{\bf\noindent Proposition 2.1.}  {\it It holds that, for all $(x,y)\in(\Omega\setminus\overline{D})^2$
\begin{equation}
%$$
\begin{array}{ll}
\displaystyle
\nabla\cdot\mbox{\boldmath $w$}_x(y)+\nabla\cdot\mbox{\boldmath $w$}_x^1(y)
&
\displaystyle
={\mathcal I}(x,y)+{\mathcal I}^1(x,y)
\\
\\
\displaystyle
&
\displaystyle
\,\,\,+\int_{\partial\Omega}\nabla G(z-x)\cdot\frac{\partial}{\partial\nu}\mbox{\boldmath $w$}_y
-\nabla G(z-y)\cdot\frac{\partial}{\partial\nu}\mbox{\boldmath $w$}_x\,dS(z),
\end{array}
\tag {2.1}
%$$
\end{equation}
where
\begin{equation}
%$$
\displaystyle
{\mathcal I}(x,y)=\int_{\Omega}V(z)\mbox{\boldmath $w$}_x\cdot\nabla G(z-y)+V(z)\nabla G(z-x)\cdot\nabla G(z-y)\,dz
\tag {2.2}
%$$
\end{equation}
and
\begin{equation}
%$$
\displaystyle
{\mathcal I}^1(x,y)=-\int_{\Omega}\nabla\mbox{\boldmath $w$}_x^1\cdot\nabla\mbox{\boldmath $w$}_y^1
+V(z)\mbox{\boldmath $w$}_x^1\cdot\mbox{\boldmath $w$}_y^1\,dz
+\int_{\partial\Omega}\nabla G(z-x)\cdot\frac{\partial}{\partial\nu}\nabla G(z-y)\,dS(z).
\tag {2.3}
%$$
\end{equation}
}

$\quad$

\noindent
From (2.2) we see ${\mathcal I}(x,y)$ at $x=y$ coincides with ${\mathcal I}(x)$ defined by (1.4)
; from (2.3)  ${\mathcal I}^1(x,y)$ at $x=y$ coincides with ${\mathcal I}^1(x)$ defined by (1.8).
Thus taking the value of the both-side on (2.1) at $x=y$,  we obtain the decomposition formula (1.6).
This completes the Proof of Theorem 1.2.

Thus remaining and main task of this subsection is to give a proof of Proposition 2.1.

$\quad$

$\quad$

{\it\noindent Proof of Proposition 2.1.}
Fix $x\in\Omega\setminus\overline{D}$.
First by elliptic regularity theory \cite{Gr}, we have 
$\mbox{\boldmath $w$}_x^1\in H^2(\Omega\setminus\overline{D})^3$ under regularity assumption $\partial\Omega\in C^{1,1}$.
Besides, by assumption (1.2) the function $\mbox{\boldmath $w$}_x^1$ is harmonic in $\Omega\setminus\overline{D}$
and thus smooth in a neighbourhood of $y\in\Omega\setminus\overline{D}$.
Then by Green's theorem \cite{Gr} we obtain, for all $y\in\Omega\setminus\overline{D}$
\begin{equation}
%$$
\begin{array}{ll}
\displaystyle
\mbox{\boldmath $w$}_x^1(y)
&
\displaystyle
=
\int_{\partial\Omega}\frac{\partial}{\partial\nu}\mbox{\boldmath $w$}_x^1(z)G(z-y)
-\nabla G(z-x)\frac{\partial}{\partial\nu}G(z-y)\,dS(z)
\\
\\
\displaystyle
&
\displaystyle
\,\,\,
-\int_{\Omega}\mbox{\boldmath $w$}_x^1(z)V(z)G(z-y)\,dz.
\end{array}
\tag {2.4}
%$$
\end{equation}

\noindent
Here,  we go beyond the proof of Theorem 3 presented in Appendix A  of \cite{IPS}.
Differentiating both sides on (2.4) with respect to $y_j$ for $j=1,2,3$, one gets
\begin{equation}
%$$
\begin{array}{ll}
\displaystyle
\frac{\partial}{\partial y_j}\mbox{\boldmath $w$}_x^1(y)
&
\displaystyle
=
-\int_{\partial\Omega}\frac{\partial}{\partial\nu}\mbox{\boldmath $w$}_x^1(z)\frac{\partial G}{\partial z_j}(z-y)\,dS(z)
+\int_{\partial\Omega}\nabla G(z-x)\frac{\partial}{\partial\nu}\frac{\partial G}{\partial z_j}(z-y)\,dS(z)
\\
\\
\displaystyle
&
\displaystyle
\,\,\,
+\int_{\Omega}\mbox{\boldmath $w$}_x^1(z)V(z)\frac{\partial G}{\partial z_j}(z-y)\,dz.
\end{array}
\tag {2.5}
%$$
\end{equation}
Note that we made use of assumption (1.2) and $y\in\Omega\setminus\overline{D}$ in differentiating the last term of right-hand side on (2.1) with respect to $y_j$.
In what follows, we denote by $\mbox{\boldmath $A$}_i$ with $i=1,2,3$ the $i$-th component of the vector valued function $\mbox{\boldmath $A$}$.

Let us focus on the first and third term of the right-hand side on (2.5).
Let $i=1,2,3$.
Using integration by parts, from (1.7) one can write the first term as
\begin{equation}
%$$
\begin{array}{ll}
\displaystyle
-\int_{\partial\Omega}\frac{\partial}{\partial\nu}(\mbox{\boldmath $w$}_x^1)_i(z)\frac{\partial G}{\partial z_j}(z-y)\,dS(z)
&
\displaystyle
=-\int_{\partial\Omega}\frac{\partial}{\partial\nu}(\mbox{\boldmath $w$}_x^1)_i(z)(\mbox{\boldmath $w$}_y^1)_j\,dS(z)
\\
\\
\displaystyle
&
\displaystyle
=-\left(\int_{\Omega}\Delta(\mbox{\boldmath $w$}_x^1)_i(\mbox{\boldmath $w$}_y^1)_j\,dz
+
\int_{\Omega}\nabla(\mbox{\boldmath $w$}_x^1)_i
\cdot\nabla(\mbox{\boldmath $w$}_y^1)_j\,dz\right)
\\
\\
\displaystyle
&
\displaystyle
=-\int_{\Omega}(\mbox{\boldmath $w$}_x^1)_iV(z)(\mbox{\boldmath $w$}_y^1)_j\,dz
-\int_{\Omega}\nabla(\mbox{\boldmath $w$}_x^1)_i
\cdot\nabla(\mbox{\boldmath $w$}_y^1)_j\,dz.
\end{array}
\tag {2.6}
%$$
\end{equation}

Next, from (1.5) and (1.7) we see that the third term has the expression
\begin{equation}
%$$
\begin{array}{ll}
\displaystyle
\int_{\Omega}(\mbox{\boldmath $w$}_x^1)_i(z)V(z)\frac{\partial G}{\partial z_j}(z-y)\,dz
&
\displaystyle
=\int_{\Omega}(\mbox{\boldmath $w$}_x^1)_i(\Delta(\mbox{\boldmath $w$}_y)_j(z)-V(z)(\mbox{\boldmath $w$}_y)_j(z))\,dz
\\
\\
\displaystyle
&
\displaystyle
=\int_{\partial\Omega}(\mbox{\boldmath $w$}_x^1)_i\frac{\partial}{\partial\nu}(\mbox{\boldmath $w$}_y)_j\,dS(z)
-\int_{\Omega}\nabla(\mbox{\boldmath $w$}_x^1)_i\cdot\nabla(\mbox{\boldmath $w$}_y)_j\,dz
\\
\\
\displaystyle
&
\displaystyle
\,\,\,
-\int_{\Omega}(\mbox{\boldmath $w$}_x^1)_iV(z)(\mbox{\boldmath $w$}_y)_j\,dz
\\
\\
\displaystyle
&
\displaystyle
=\int_{\partial\Omega}\frac{\partial G}{\partial z_i}(z-x)\frac{\partial}{\partial\nu}(\mbox{\boldmath $w$}_y)_j\,dS(z)
+\int_{\Omega}\Delta(\mbox{\boldmath $w$}_x^1)_{i}(\mbox{\boldmath $w$}_y)_j\,dz
\\
\\
\displaystyle
&
\displaystyle
\,\,\,
-\int_{\Omega}(\mbox{\boldmath $w$}_x^1)_iV(z)(\mbox{\boldmath $w$}_y)_j\,dz
\\
\\
\displaystyle
&
\displaystyle
=\int_{\partial\Omega}\frac{\partial G}{\partial z_i}(z-x)\frac{\partial}{\partial\nu}(\mbox{\boldmath $w$}_y)_j\,dS(z).
\end{array}
\tag {2.7}
%$$
\end{equation}
Substituting (2.6) and (2.7) into (2.5), we obtain
\begin{equation}
%$$
\begin{array}{ll}
\displaystyle
\frac{\partial}{\partial y_j}(\mbox{\boldmath $w$}_x^1)_i(y)
&
\displaystyle
=-\int_{\Omega}\nabla(\mbox{\boldmath $w$}_x^1)_i\cdot\nabla(\mbox{\boldmath $w$}_y^1)_j\,dz
-\int_{\Omega}(\mbox{\boldmath $w$}_x^1)_iV(z)(\mbox{\boldmath $w$}_y^1)_j\,dz
\\
\\
\displaystyle
&
\displaystyle
\,\,\,
+\int_{\partial\Omega}\frac{\partial G}{\partial z_i}(z-x)\frac{\partial}{\partial\nu}\frac{\partial G}{\partial z_j}(z-y)\,dS(z)
+
\int_{\partial\Omega}\frac{\partial G}{\partial z_i}(z-x)\frac{\partial}{\partial\nu}(\mbox{\boldmath $w$}_y)_j\,dS(z).
\end{array}
\tag {2.8}
%$$
\end{equation}
Letting $i=j$ of the both sides on (2.8) and taking the summation with respect to $i$, we obtain
\begin{equation}
%$$
\begin{array}{ll}
\displaystyle
\nabla\cdot\mbox{\boldmath $w$}_x^1(y)
&
\displaystyle
=-\int_{\Omega}\nabla\mbox{\boldmath $w$}_x^1\cdot\nabla\mbox{\boldmath $w$}_y^1+V(z)\mbox{\boldmath $w$}_x^1\cdot\mbox{\boldmath $w$}_y^1\,dz
\\
\\
\displaystyle
&
\displaystyle
\,\,\,
+\int_{\partial\Omega}\nabla G(z-x)\cdot\frac{\partial}{\partial\nu}\nabla G(z-y)\,dS(z)
+\int_{\partial\Omega}\,\nabla G(z-x)\cdot\frac{\partial}{\partial\nu}\mbox{\boldmath $w$}_y\,dS(z).
\end{array}
\tag {2.9}
%$$
\end{equation}

Similarly to (2.4), we have also the expression for $\mbox{\boldmath $w$}_x$:
\begin{equation}
%$$
\begin{array}{ll}
\displaystyle
\mbox{\boldmath $w$}_x(y)
&
\displaystyle
=
\int_{\partial\Omega}\frac{\partial}{\partial\nu}\mbox{\boldmath $w$}_x(z)G(z-y)\,dS(z)
\\
\\
\displaystyle
&
\displaystyle
\,\,\,
-\int_{\Omega}\mbox{\boldmath $w$}_x(z)V(z)G(z-y)\,dz
-\int_{\Omega}V(z)\nabla G(z-x)G(z-y)\,dz.
\end{array}
\tag {2.10}
%$$
\end{equation}
Differentiating both sides of equation on (2.10) with respect to $y_j$, we obtain, for $i=1,2,3$
\begin{equation}
%$$
\begin{array}{ll}
\displaystyle
\frac{\partial}{\partial y_j}(\mbox{\boldmath $w$}_x)_i(y)
&
\displaystyle
=
-\int_{\partial\Omega}\frac{\partial}{\partial\nu}(\mbox{\boldmath $w$}_x)_i(z)\frac{\partial G}{\partial z_j}(z-y)\,dS(z)
\\
\\
\displaystyle
&
\displaystyle
\,\,\,
+\int_{\Omega}(\mbox{\boldmath $w$}_x)_i(z)V(z)\frac{\partial G}{\partial z_j}(z-y)\,dz
+\int_{\Omega}V(z)\frac{\partial G}{\partial z_i}(z-x)\frac{\partial G}{\partial z_j}(z-y)\,dz.
\end{array}
\tag {2.11}
%$$
\end{equation}
Letting $i=j$ of the both sides on (2.11) and taking the summation with respect to $i$, we obtain
\begin{equation}
%$$
\begin{array}{ll}
\displaystyle
\nabla\cdot\mbox{\boldmath $w$}_x(y)
&
\displaystyle
=-\int_{\partial\Omega}\frac{\partial}{\partial\nu}\mbox{\boldmath $w$}_x\cdot\nabla G(z-y)\,dS(z)
\\
\\
\displaystyle
&
\displaystyle
\,\,\,
+\int_{\Omega}V(z)\mbox{\boldmath $w$}_x\cdot\nabla G(z-y)+V(z)\nabla G(z-x)\cdot\nabla G(z-y)\,dz.
\end{array}
\tag {2.12}
%$$
\end{equation}
Now introducing ${\mathcal I}(x,y)$ and ${\mathcal I}^1(x,y)$ by (2.2) and (2.3), respectively, from (2.9) and (2.12) one obtains (2.1).

\noindent
$\Box$

$\quad$

{\bf\noindent Remark 2.1.}
The original method of proof of (1.18) of Theorem 3 in \cite{IPS} which corresponds to (1.6)
employs the function $G(\,\cdot\,-x)-v_n$, which is still a fundamental solution of the Laplace equation in the sense of distribution
and satisfies $G(\,\cdot\,-x)-v_n\rightarrow 0$ in $H^1(B)$ for all open ball $B$ with $\overline{B}\subset\Omega\setminus\sigma$.
This means the distribution $G(\,\cdot\,-x)-v_n$ plays a role of the Carleman function, which has been pointed out in \cite{ICar}.
However, the function $G_n^j(\,\cdot\,,x)$ defined by (1.12) is not a fundamental solution as a distribution.
So we follow a direct proof of (1.18) of Theorem 3, as presented in Appendix A of \cite{IPS}, which neither makes use of this convergence property
nor involves the sequence $\{v_n^j\}$.
This approach has been employed also in \cite{IPS2} and \cite{IPS3}.

\subsection{Remarks on Liftings and Proof of Proposition 1.1}

$\quad$

We call ${\mathcal I}(x,y)$ and ${\mathcal I}^1(x,y)$  defined by (2.2) and (2.3) the {\it liftings} of ${\mathcal I}(x)$ and ${\mathcal I}^1(x)$, respectively.
In this subsection we study the properties of lifings of ${\mathcal I}(x)$ and ${\mathcal I}^1(x)$ together with $\nabla\cdot\mbox{\boldmath $w$}_x(y)$ and $\nabla\cdot\mbox{\boldmath $w$}_x^1(y)$.

$\quad$

{\bf\noindent On ${\mathcal I}^1(x,y)$ and $\mbox{\boldmath $w$}_x^1(y)$.}

$\quad$

Let $(x,y)\in\,(\Omega\setminus\overline{D})^2$.
Green's theorem outside $\Omega$ yields
$$\begin{array}{ll}
\displaystyle
-\int_{\partial\Omega}\nabla G(z-x)\cdot\frac{\partial}{\partial\nu}\nabla G(z-y)\,dS(z)
&
\displaystyle
=\int_{\Bbb R^3\setminus\overline{\Omega}}\nabla G(z-x)\cdot\Delta\nabla G(z-y)\,dz
\\
\\
\displaystyle
&
\displaystyle
\,\,\,
+\int_{\Bbb R^3\setminus\overline{\Omega}}\nabla^2G(z-x)\cdot\nabla^2G(z-y)\,dz
\\
\\
&
\displaystyle
=\int_{\Bbb R^3\setminus\overline{\Omega}}\nabla^2G(z-x)\cdot\nabla^2G(z-y)\,dz.
\end{array}
$$
Thus (2.3) has also an energy integral expression
\begin{equation}
%$$
\displaystyle
{\mathcal I}^1(x,y)=-\int_{\Omega}\nabla\mbox{\boldmath $w$}_x^1\cdot\nabla\mbox{\boldmath $w$}_y^1+V(z)\mbox{\boldmath $w$}_x^1\cdot\mbox{\boldmath $w$}_y^1\,dz
-\int_{\Bbb R^3\setminus\overline{\Omega}}\nabla^2 G(z-x)\cdot\nabla^2 G(z-y)\,dz.
\tag {2.13}
%$$
\end{equation}
This yields the symmetry
\begin{equation}
%$$
\displaystyle
{\mathcal I}^1(x,y)={\mathcal I}^1(y,x).
\tag {2.14}
%$$
\end{equation}
And we have the symmetry
\begin{equation}
%$$
\displaystyle
\int_{\partial\Omega}\nabla G(z-x)\cdot\frac{\partial}{\partial\nu}\nabla G(z-y)\,dS(z)
=\int_{\partial\Omega}\nabla G(z-y)\cdot\frac{\partial}{\partial\nu}\nabla G(z-x)\,dS(z).
\tag {2.15}
%$$
\end{equation}

One can calculate ${\mathcal I}^1(x,y)$ from $\Lambda_V$ directly as follows.   
From equation (2.6) we have
\begin{equation}
%$$
\displaystyle
\int_{\Omega}\nabla\mbox{\boldmath $w$}_x^1\cdot\nabla\mbox{\boldmath $w$}_y^1+V(z)\mbox{\boldmath $w$}_x^1\cdot\mbox{\boldmath $w$}_y^1\,dz
=\int_{\partial\Omega}\frac{\partial}{\partial\nu}\mbox{\boldmath $w$}_x^1\cdot\nabla G(z-y)\,dS(z).
\tag {2.16}
%$$
\end{equation}
By the definition of the Dirichlet-to-Neumann map, one has the expression:
\begin{equation}
%$$
\displaystyle
\left.\frac{\partial}{\partial\nu}\mbox{\boldmath $w$}_x^1\right\vert_{\partial\Omega}=\left.\Lambda_V\left(\nabla G(\,\cdot\,-x)\right\vert_{\partial\Omega}\right).
\tag {2.17}
%$$
\end{equation}
Thus, from (2.3), (2.15), (2.16) and (2.17) we obtain the calculation formula of ${\mathcal I}^1(x,y)$ from $\Lambda_V$:
\begin{equation}
%$$
\displaystyle
{\mathcal I}^1(x,y)=
\int_{\partial\Omega}\left\{\frac{\partial}{\partial\nu}\nabla G(z-x)-\Lambda_V(\nabla G(\,\cdot\,-x)\vert_{\partial\Omega})\right\}\cdot\nabla G(z-y)\,dS(z).
\tag {2.18}
%$$
\end{equation}
Besides, formula (2.18) enables us to conclude that, for each fixed $x\in\Omega\setminus\overline{D}$
\begin{equation}
%$$
\begin{array}{ll}
\displaystyle
\Delta_y{\mathcal I}^1(x,y)=0, & y\in\Omega\setminus\overline{D}
\end{array}
\tag {2.19}
%$$
\end{equation}
and (2.14) yields, for each fixed $y\in\Omega\setminus\overline{D}$
\begin{equation}
%$$
\begin{array}{ll}
\displaystyle
\Delta_x{\mathcal I}^1(x,y)=0, & x\in\Omega\setminus\overline{D}.
\end{array}
\tag {2.20}
%$$
\end{equation}
Note that $\Delta_y$ and $\Delta_x$ denote the Laplacian with respect to independent variables $y$ and $x$, respectively.

From (2.9) and (2.3) we have another expression
\begin{equation}
%$$
\displaystyle
{\mathcal I}^1(x,y)=\nabla\cdot\mbox{\boldmath $w$}_x^1(y)-
\int_{\partial\Omega}\nabla G(z-x)\cdot\frac{\partial}{\partial\nu}\mbox{\boldmath $w$}_y\,dS(z).
\tag {2.21}
%$$
\end{equation}
This yields
$$\displaystyle
\int_{\partial\Omega}\nabla G(z-x)\cdot\frac{\partial}{\partial\nu}\mbox{\boldmath $w$}_y\,dS(z)
=\nabla\cdot\mbox{\boldmath $w$}_x^1(y)-{\mathcal I}^1(x,y).
$$
Then, from (2.19) and assumption (1.2)  one gets, for each fixed $x\in\Omega\setminus\overline{D}$
\begin{equation}
%$$
\begin{array}{ll}
\displaystyle
\Delta_y\int_{\partial\Omega}\nabla G(z-x)\cdot\frac{\partial}{\partial\nu}\mbox{\boldmath $w$}_y\,dS(z)=0, & y\in\Omega\setminus\overline{D}.
\end{array}
\tag {2.22}
%$$
\end{equation}
It follows from (2.21) and (2.14) one gets
$$\displaystyle
\nabla\cdot\mbox{\boldmath $w$}_y^1(x)
=\nabla\cdot\mbox{\boldmath $w$}_x^1(y)+\int_{\partial\Omega}\nabla G(z-y)\cdot\frac{\partial}{\partial\nu}\mbox{\boldmath $w$}_x
-\nabla G(z-x)\frac{\partial}{\partial\nu}\mbox{\boldmath $w$}_y\,dS(z).
$$
This together with (2.22) and (1.2) yields, for each fixed $x\in\Omega\setminus\overline{D}$
$$\begin{array}{ll}
\displaystyle
\Delta_y\nabla\cdot\mbox{\boldmath $w$}_y^1(x)=0, & y\in\Omega\setminus\overline{D}
\end{array}
$$
and thus, for each fixed $y\in\Omega\setminus\overline{D}$
\begin{equation}
%$$
\begin{array}{ll}
\displaystyle
\Delta_x\nabla\cdot\mbox{\boldmath $w$}_x^1(y)=0, & x\in\Omega\setminus\overline{D}.
\end{array}
\tag {2.23}
%$$
\end{equation}

$\quad$

{\bf\noindent On ${\mathcal I}(x,y)$ and $\nabla\cdot\mbox{\boldmath $w$}_x(y)$.}

$\quad$

Let $(x,y)\in\,(\Omega\setminus\overline{D})^2$ and consider the symmetry of ${\mathcal I}(x,y)$.
Using integration by parts and (1.5), we have
$$\begin{array}{ll}
\displaystyle
-\int_{\Omega}V(z)\mbox{\boldmath $w$}_x\cdot\nabla G(z-y)\,dz
&
\displaystyle
=-\int_{\Omega}\mbox{\boldmath $w$}_x\cdot V(z)\nabla G(z-y)\,dz
\\
\\
\displaystyle
&
\displaystyle
=\int_{\Omega}\mbox{\boldmath $w$}_x\cdot(-\Delta+V(z))\mbox{\boldmath $w$}_y\,dz
\\
\\
\displaystyle
&
\displaystyle
=-\int_{\partial\Omega}\mbox{\boldmath $w$}_x\cdot\frac{\partial}{\partial\nu}\mbox{\boldmath $w$}_y\,dS(z)
+\int_{\Omega}\nabla\mbox{\boldmath $w$}_x\cdot\nabla\mbox{\boldmath $w$}_y+V(z)\mbox{\boldmath $w$}_x\cdot\mbox{\boldmath $w$}_y\,dz
\\
\\
\displaystyle
&
\displaystyle
=\int_{\Omega}\nabla\mbox{\boldmath $w$}_x\cdot\nabla\mbox{\boldmath $w$}_y+V(z)\mbox{\boldmath $w$}_x\cdot\mbox{\boldmath $w$}._y\,dz.
\end{array}
$$
From this together with (2.2) we obtain
\begin{equation}
%$$
\displaystyle
{\mathcal I}(x,y)=-\int_{\Omega}\nabla\mbox{\boldmath $w$}_x\cdot\nabla\mbox{\boldmath $w$}_y+V(z)\mbox{\boldmath $w$}_x\cdot\mbox{\boldmath $w$}_y\,dz+\int_{\Omega}V(z)\nabla G(z-x)\cdot\nabla G(z-y)\,dz.
\tag {2.24}
%$$
\end{equation}
Thus one gets the symmetry
\begin{equation}
%$$
\begin{array}{ll}
\displaystyle
{\mathcal I}(x,y)={\mathcal I}(y,x), & (x,y)\in(\Omega\setminus\overline{D})^2.
\end{array}
\tag {2.25}
%$$
\end{equation}
Besides, it follows from (2.12) and (2.2) that
\begin{equation}
%$$
\displaystyle
{\mathcal I}(x,y)=\nabla\cdot\mbox{\boldmath $w$}_x(y)+\int_{\partial\Omega}\frac{\partial}{\partial\nu}\mbox{\boldmath $w$}_x\cdot\nabla G(z-y)\,dS(z).
\tag {2.26}
%$$
\end{equation}
This together with assumption (1.2), we see, for each fixed $x\in\Omega\setminus\overline{D}$,
\begin{equation}
%$$
\begin{array}{ll}
\displaystyle
\Delta_y{\mathcal I}(x,y)=0, & y\in\Omega\setminus\overline{D}
\end{array}
\tag {2.27}
%$$
\end{equation}
and symmetry (2.25) yields, for each fixed $y\in\Omega\setminus\overline{D}$,
\begin{equation}
%$$
\begin{array}{ll}
\displaystyle
\Delta_x{\mathcal I}(x,y)=0, & x\in\Omega\setminus\overline{D}.
\end{array}
\tag {2.28}
%$$
\end{equation}
Changing the role of $x$ and $y$ in (2.22), we have,
$$\begin{array}{ll}
\displaystyle
\Delta_x\int_{\partial\Omega}\nabla G(z-y)\cdot\frac{\partial}{\partial\nu}\mbox{\boldmath $w$}_x\,dS(z)=0, & x\in\Omega\setminus\overline{D}.
\end{array}
$$
This together with (2.26) and (2.28) yields, for each fixed $y\in\Omega\setminus\overline{D}$
\begin{equation}
%$$
\begin{array}{ll}
\Delta_x\nabla\cdot\mbox{\boldmath $w$}_x(y)=0, & x\in\Omega\setminus\overline{D}.
\end{array}
\tag {2.29}
%$$
\end{equation}

\subsubsection{Proof of Proposition 1.1}

$\quad$

Let $(x,y)\in(\Omega\setminus\overline{D})^2$.
Integration by parts yields
$$\begin{array}{ll}
\displaystyle
-\int_{\Omega}\nabla\mbox{\boldmath $w$}_x\cdot\nabla\mbox{\boldmath $w$}_y^1\,dz
&
\displaystyle
=\int_{\Omega}
\mbox{\boldmath $w$}_x\cdot\Delta\mbox{\boldmath $w$}_y^1\,dz
\\
\\
\displaystyle
&
\displaystyle
=\int_{\Omega}\mbox{\boldmath $w$}_x\cdot V(z)\mbox{\boldmath $w$}_y^1\,dz.
\end{array}
$$
Thus one gets the orthogonality relation between $\mbox{\boldmath $w$}_x$ and $\mbox{\boldmath $w$}_y^1$ in the following sense:
$$\displaystyle
\int_{\Omega}\nabla\mbox{\boldmath $w$}_x\cdot\nabla\mbox{\boldmath $w$}_y^1+V(z)\mbox{\boldmath $w$}_x\cdot\mbox{\boldmath $w$}_y^1\,dz=0.
$$
This implies
$$\begin{array}{l}
\displaystyle
\,\,\,\,\,\,
\int_{\Omega}\nabla\mbox{\boldmath $W$}_x\cdot\nabla\mbox{\boldmath $W$}_y+V(z)\mbox{\boldmath $W$}_x\cdot\mbox{\boldmath $W$}_y\,dz
\\
\\
\displaystyle
=\int_{\Omega}\nabla\mbox{\boldmath $w$}_x\cdot\nabla\mbox{\boldmath $w$}_y+V(z)\mbox{\boldmath $w$}_x\cdot\mbox{\boldmath $w$}_y\,dz
+\int_{\Omega}\nabla\mbox{\boldmath $w$}_x^1\cdot\nabla\mbox{\boldmath $w$}_y^1+V(z)\mbox{\boldmath $w$}_x^1\cdot\mbox{\boldmath $w$}_y^1\,dz.
\end{array}
$$
Therefore from this, (2.13) and (2.24) we obtain
$$\begin{array}{l}
\displaystyle
\,\,\,\,\,\,
-\int_{\Omega}\nabla\mbox{\boldmath $W$}_x\cdot\nabla\mbox{\boldmath $W$}_y+V(z)\mbox{\boldmath $W$}_x\cdot\mbox{\boldmath $W$}_y\,dz
\\
\\
\displaystyle
={\mathcal I}^1(x,y)+\int_{\Bbb R^3\setminus\overline{\Omega}}\nabla^2 G(z-x)\cdot\nabla^2 G(z-y)\,dz
\\
\\
\displaystyle
\,\,\,
+{\mathcal I}(x,y)-\int_{\Omega}V(z)\nabla G(z-x)\cdot\nabla G(z-y)\,dz.
\end{array}
$$
Then, by (2.1) and the definition of $\mbox{\boldmath $W$}_x(y)$, we obtain
$$\begin{array}{ll}
\displaystyle
\nabla\cdot\mbox{\boldmath $W$}_x(y)
&
\displaystyle
=-\int_{\Omega}\nabla\mbox{\boldmath $W$}_x\cdot\nabla\mbox{\boldmath $W$}_y+V(z)\mbox{\boldmath $W$}_x\cdot\mbox{\boldmath $W$}_y\,dz
\\
\\
\displaystyle
&
\displaystyle
\,\,\,
+\int_{\Omega} V(z)\nabla G(z-x)\cdot\nabla G(z-y)\,dz-\int_{\Bbb R^3\setminus\overline{\Omega}}\nabla^2 G(z-x)\cdot\nabla^2 G(z-y)\,dz
\\
\\
\displaystyle
&
\displaystyle
+\int_{\partial\Omega}\nabla G(z-x)\cdot\frac{\partial}{\partial\nu}\mbox{\boldmath $w$}_y
-\nabla G(z-y)\cdot\frac{\partial}{\partial\nu}\mbox{\boldmath $w$}_x\,dS(z).
\end{array}
$$
Taking the symmetric part of both sides, we complete the proof of Proposition 1.1.

\noindent
$\Box$

\subsection{Proof of Theorem 1.3}

$\quad$

{\it\noindent Proof of {\rm(a)}}.
Given $(x,y)\in(\Omega\setminus\overline{D})^2$ and $\sigma\in N_x$, $\sigma'\in N_y$ satisfying $\sigma\cap\overline{D}=\sigma'\cap\overline{D}=\emptyset$, 
let $\{v_n^j\}$ and $\{v_n'^j\}$
be needle sequences for $(x,\sigma)$ and $(y,\sigma')$, respectively.  
By the expression (2.2), similar to the derivation of (1.3) in \cite{IProbe}, we know
\begin{equation}
%$$
\displaystyle
{\mathcal I}(x,y)=\lim_{n\rightarrow\infty}\sum_{j}\int_{\partial\Omega}(\Lambda_V-\Lambda_0)v_n^j\vert_{\partial\Omega}\cdot v_n'^j\vert_{\partial\Omega}\,dS(z)
\tag {2.30}
%$$
\end{equation}
and also
\begin{equation}
%$$
\displaystyle
\int_{\partial\Omega}\frac{\partial}{\partial\nu}\mbox{\boldmath $w$}_x\cdot\nabla G(z-y)\,dS(z)
=\lim_{n\rightarrow\infty}\sum_{j}\int_{\partial\Omega}(\Lambda_V-\Lambda_0)v_n^j\vert_{\partial\Omega}\cdot\frac{\partial G}{\partial z_j}(z-y)\,dS(z).
\tag {2.31}
%$$
\end{equation}
Then, from (2.26) we obtain
\begin{equation}
%$$
\begin{array}{ll}
\displaystyle
\nabla\cdot\mbox{\boldmath $w$}_x(y)
&
\displaystyle
=\lim_{n\rightarrow\infty}\sum_j
\int_{\partial\Omega}
(\Lambda_V-\Lambda_0)v_n^j\vert_{\partial\Omega}\cdot (v_n'^j(z)-\frac{\partial G}{\partial z_j}(z-y))\,dS(z)
\\
\\
\displaystyle
&
\displaystyle
=-\lim_{n\rightarrow\infty}\sum_j
\int_{\partial\Omega}
(\Lambda_V-\Lambda_0)v_n^j\vert_{\partial\Omega}\cdot (\frac{\partial G}{\partial z_j}(z-y)-v_n'^j(z))\,dS(z).
\end{array}
\tag {2.32}
%$$
\end{equation}
Setting $y=x$ and choosing such as $v_n'^j=v_n^j$ in (2.32), one gets (1.13).

$\quad$

{\it\noindent Proof of {\rm(b)}}.
Estimates (1.14) and (1.15) are proved as follows.
It follows from (2.5) and (2.17) together with regularity assumption $\partial\Omega\in C^{1,1}$
that
\begin{equation}
%$$
\begin{array}{ll}
\displaystyle
\nabla\cdot\mbox{\boldmath $w$}_x^1(x)
&
\displaystyle
=-<\Lambda_V(\nabla G(\,\cdot\,-x)\vert_{\partial\Omega},\nabla G(\,\cdot\,-x)\vert_{\partial\Omega}>
\\
\\
\displaystyle
&
\displaystyle
\,\,\,
+\int_{\partial\Omega}\nabla G(z-x)\cdot\frac{\partial}{\partial\nu}\nabla G(z-x)\,dS(z)
+\int_{\Omega}\mbox{\boldmath $w$}_x^1(z)\cdot V(z)\nabla G(z-x)\,dz.
\\
\\
\displaystyle
&
\displaystyle
\equiv {\rm I}+{\rm II}+{\rm III}.
\end{array}
\tag {2.33}
%$$
\end{equation}
Here we know
$$\displaystyle
\sup_{x\in\Omega\setminus\overline{D},\,\text{dist}\,(x,\partial\Omega)>\epsilon}\,\Vert\nabla G(\,\cdot\,-x)\vert_{\partial\Omega}\Vert_{H^{\frac{1}{2}}(\partial\Omega)}<\infty.
$$
Thus 
\begin{equation}
%$$
\displaystyle
\sup_{x\in\Omega\setminus\overline{D},\,\text{dist}\,(x,\partial\Omega)>\epsilon}
\vert{ \rm I} \vert <\infty,
\tag {2.34}
%$$
\end{equation}
Besides by estimating the integrand point wisely, one has
\begin{equation}
%$$
\displaystyle
\sup_{x\in\Omega\setminus\overline{D},\, \text{dist}\,(\partial\Omega,x)>\epsilon}\vert{\rm II}\vert<\infty.
\tag {2.35}
%$$
\end{equation}

For estimating the third term ${\rm III}$,  we know
$$\displaystyle
\sup_{x\in\Omega\setminus\overline{D}\,, \text{dist}\,(\partial\Omega,x)>\epsilon}\,
\Vert\nabla G(\,\cdot\,-x)\Vert_{\partial\Omega}\Vert_{H^{\frac{3}{2}}(\partial\Omega)}<\infty.
$$
Thus, by the well-posedness of the Dirichlet problem (1.1) in $H^2(\Omega)$ one gets
\begin{equation}
%$$
\displaystyle
\sup_{x\in\Omega\setminus\overline{D},\,\text{dist}\,(x,\partial\Omega)>\epsilon}\,\Vert\mbox{\boldmath $w$}_x^1\Vert_{H^2(\Omega)}<\infty.
\tag {2.36}
%$$
\end{equation}
Then the Sobolev embedding yields
$$\displaystyle
C=\sup_{x\in\Omega\setminus\overline{D},\,\text{dist}\,(x,\partial\Omega)>\epsilon}\,\Vert\mbox{\boldmath $w$}_x^1\Vert_{L^{\infty}(\Omega)}<\infty.
$$
This together with (1.2), one gets
$$\displaystyle
\vert{\rm III}\vert
\le C\Vert V\Vert_{L^{\infty}(D)}\Vert\nabla G(\,\cdot\,-x)\Vert_{L^1(D)}.
$$
Since
$$\displaystyle
\sup_{x\in\Bbb R^3}\,\Vert\nabla G(\,\cdot\,-x)\Vert_{L^1(D)}<\infty,
$$
we conclude
\begin{equation}
%$$
\displaystyle
\sup_{x\in\Omega\setminus\overline{D},\,\text{dist}\,(x,\partial\Omega)>\epsilon}\vert{\rm III}\vert
<\infty.
\tag {2.37}
%$$
\end{equation}
Therefore from (2.33), (2.34), (2.35) and (2.37) we obtain (1.14).  Estimate (1.15) is a consequence of expression (1.8) and estimates (2.35), (2.36).
Finally (1.16) follows from the decomposition (1.6) combined with estimates  (1.14) and (1.15).

$\quad$

{\it\noindent Proof of {\rm(c)}}.
The proof of {\rm (c)} is a direct consequence of {\rm(b)} of Theorem 1.1 and  the estimate (1.16).

$\quad$

{\it\noindent Proof of {\rm(d)}}.  This is a consequence of (1.16) and {\rm (c)} of Theorem 1.1.

\noindent
$\Box$

\subsection{Proof of Theorem 1.4}

$\quad$

{\it\noindent Proof of {\rm(a)}.}
The derivation of (1.19) are similar to that of (1.3).  
The key point is: the solution $u=u_n^j$ of (1.1) with $f=(v_n^j+H^j(\,\cdot\,,x))\vert_{\partial\Omega}$ has the expression
$$\displaystyle
u_n^j=v_n^j+H^j(\,\cdot\,,x)+(w^*)_n^j,
$$
where $(w^*)_n^j=w$ solves
$$\left\{
\begin{array}{ll}
\displaystyle
-\Delta w+V(z)w=-V(z)(v_n^j+H^j(\,\cdot\,,x)), & z\in\Omega,
\\
\\
\displaystyle
w=0, & z\in\partial\Omega.
\end{array}
\right.,
$$
Then, by the well posedness, we see that, as $n\rightarrow\infty$  the sequence $\{(w^*)_n^j\}$ converges to the $j$-th component of $\mbox{\boldmath $w$}^*_x$ in $H^2(\Omega)$
which is the unique solution of (1.21). 

\noindent
Besides, the Alessandrini identity together with (1.18) yields
$$\begin{array}{l}
\displaystyle
\,\,\,\,\,\,
\sum_{j}\,<(\Lambda_V-\Lambda_0)G_n^j(\,\cdot\,,x)\vert_{\partial\Omega},G_n^j(\,\cdot\,,x)\vert_{\partial\Omega}>
\\
\\
\displaystyle
=\sum_{j}\,<(\Lambda_V-\Lambda_0)(v_n^j+H^j(\,\cdot\,,x))\vert_{\partial\Omega},(v_n^j+H^j(\,\cdot\,,x))\vert_{\partial\Omega}>
\\
\\
\displaystyle
=\sum_j\,\int_D V(z)u_n^j (v_n^j+H^j(\,\cdot\,,x))\,dz
\\
\\
\displaystyle
\rightarrow \int_D V(z)(\mbox{\boldmath $G$}(z,x)+
\mbox{\boldmath $w$}^*_x(z))\cdot\mbox{\boldmath $G$}(z,x)\,dz,
\end{array}
$$
where the $j$-th component of $\mbox{\boldmath $G$}$ is given by $G^j$ defined by (1.17).
Thus, using ${\mathcal I}^*(x)$ defined by (1.20), one obtains the expression (1.19).

$\quad$

{\it\noindent Proof of {\rm (b)}.}
The symmetry of $\Lambda_v$ together with (1.18) yields
$$\begin{array}{l}
\displaystyle
\,\,\,\,\,\,
\sum_j<(\Lambda_V-\Lambda_0)G_n^j(\,\cdot\,,x)\vert_{\partial\Omega},G_n^j(\,\cdot\,,x)\vert_{\partial\Omega}>
\\
\\
\displaystyle
=\sum_j<(\Lambda_V-\Lambda_0)v_n^j\vert_{\partial\Omega},v_n^j\vert_{\partial\Omega}>
\\
\\
\displaystyle
\,\,\,
+2\sum_j<(\Lambda_V-\Lambda_0)v_n^j\vert_{\partial\Omega}, H^j(\,\cdot\,,x)\vert_{\partial\Omega}>
+\sum_j<(\Lambda_V-\Lambda_0)H^j(\,\cdot\,,x)\vert_{\partial\Omega},H^j(\,\cdot\,,x)\vert_{\partial\Omega}>.
\end{array}
$$
Taking the limit of both sides, from (1.3), (2.31) with $y=x$ and the boundary condition of $\mbox{\boldmath $H$}(\,\cdot\,,x)$ on $\partial\Omega$ we conclude
$$\displaystyle
{\mathcal I}^*(x)={\mathcal I}(x)-2\int_{\partial\Omega}\frac{\partial}{\partial\nu}\mbox{\boldmath $w$}_x\cdot\nabla G(z-x)\,dS(z)
+<(\Lambda_V-\Lambda_0)\nabla G(\,\cdot\,,-x)\vert_{\partial\Omega}, \nabla G(\,\cdot\,-x)\vert_{\partial\Omega}>.
$$
Then from (2.21) we obtain (1.22).

$\quad$

{\it\noindent Proof of {\rm (c)}.}
It follows from (1.14), (1.15), (1.22) and (2.34) (also (2.34) with $V\equiv 0$) that, for each $\epsilon>0$
\begin{equation}
%$$
\displaystyle
\sup_{x\in\Omega\setminus\overline{D},\,\text{dist}\,(x,\partial\Omega)>\epsilon}\,\vert{\mathcal I}^*(x)-{\mathcal I}(x)\vert<\infty.
\tag {2.38}
%$$
\end{equation}
Thus, one gets, as $x\rightarrow a\in\partial D$
$$\displaystyle
{\mathcal I}^*(x)={\mathcal I}(x)+O(1).
$$
Therefore {\rm (c)} is equivalent to that of Theorem 1.1, and justified.

$\quad$

{\it\noindent Proof of {\rm(d)}.}
This is a direct consequence of {\rm (c)} of Theorem 1.1 combined with (2.38).

$\quad$

{\it\noindent Proof of {\rm(e)}.}
Set
$$\displaystyle
\mbox{\boldmath $R$}_x=\mbox{\boldmath $w$}^*_x-\mbox{\boldmath $w$}_x
$$
where $x\in\Omega\setminus\overline{D}$.  Then, from (1.5) and (1.21) we see that the vector valued function $\mbox{\boldmath $R$}_x=\mbox{\boldmath $R$}$ solves
$$\left\{
\begin{array}{ll}
\displaystyle
-\Delta \mbox{\boldmath $R$}+V(z)\mbox{\boldmath $R$}=-V(z)\mbox{\boldmath $H$}(z,x), & z\in\Omega,
\\
\\
\displaystyle
\mbox{\boldmath $R$}=\mbox{\boldmath $0$}, & z\in\partial\Omega.
\end{array}
\right.
$$
Green's representation yields
$$\begin{array}{ll}
\displaystyle
\mbox{\boldmath $R$}_x(y)
&
\displaystyle
=
\int_{\partial\Omega}\frac{\partial}{\partial\nu}\mbox{\boldmath $R$}_x(z)G(z-y)\,dS(z)
\\
\\
\displaystyle
&
\displaystyle
\,\,\,
-\int_{D}\mbox{\boldmath $R$}_x(z)V(z)G(z-y)\,dz
-\int_{D}V(z)\mbox{\boldmath $H$}(z,x)G(z-y)\,dz.
\end{array}
$$
Taking the divergence of the both sides for $y\in\Omega\setminus\overline{D}$, we obtain
$$\begin{array}{ll}
\displaystyle
\nabla\cdot\mbox{\boldmath $R$}_x(y)
&
\displaystyle
=-
\int_{\partial\Omega}\frac{\partial}{\partial\nu}\mbox{\boldmath $R$}_x(z)\cdot\nabla G(z-y)\,dS(z)
\\
\\
\displaystyle
&
\displaystyle
\,\,\,
+\int_{D}V(z)\mbox{\boldmath $R$}_x(z)\cdot\nabla G(z-y)\,dz
+\int_{D}V(z)\mbox{\boldmath $H$}(z,x)\cdot\nabla G(z-y)\,dz.
\end{array}
$$
Besides we have (2.12), that is
$$\begin{array}{ll}
\displaystyle
\nabla\cdot\mbox{\boldmath $w$}_x(y)
&
\displaystyle
=-\int_{\partial\Omega}\frac{\partial}{\partial\nu}\mbox{\boldmath $w$}_x\cdot\nabla G(z-y)\,dS(z)
\\
\\
\displaystyle
&
\displaystyle
\,\,\,
+\int_{D}V(z)\mbox{\boldmath $w$}_x\cdot\nabla G(z-y)\,dz
+\int_{D}V(z)\nabla G(z-x)\cdot\nabla G(z-y)\,dz.
\end{array}
$$
Thus one gets
$$\begin{array}{ll}
\displaystyle
\nabla\cdot\mbox{\boldmath $w$}^*_x(y)
&
\displaystyle
=-\int_{\partial\Omega}\frac{\partial}{\partial\nu}\mbox{\boldmath $w$}^*_x\cdot\nabla G(z-y)\,dS(z)
\\
\\
\displaystyle
&
\displaystyle
\,\,\,
+\int_{D}V(z)(\mbox{\boldmath $w$}^*_x(z)+\nabla G(z-x)+\mbox{\boldmath $H$}(z,x))\cdot\nabla G(z-y)\,dz.
\end{array}
$$
Here we have
$$\begin{array}{l}
\displaystyle
\,\,\,\,\,\,
\lim_{n\rightarrow\infty}\sum_j\,<(\Lambda_V-\Lambda_0)(v_n^j+H^j(\,\cdot\,,x))\vert_{\partial\Omega}, H^j(\,\cdot\,,x)\vert_{\partial\Omega}>
\\
\\
\displaystyle
=-\lim_{n\rightarrow\infty}\int_{\partial\Omega}\frac{\partial}{\partial\nu}(u_n^j-(v_n^j+H^j(\,\cdot\,,x))\cdot\nabla G(z-x)\,dS(z)
\\
\\
\displaystyle
=-\lim_{n\rightarrow\infty}\,\sum_j\int_{\partial\Omega}\frac{\partial}{\partial\nu}(w^*)_n^j\frac{\partial G}{\partial z_j}(z-x)\,dS(z)
\\
\\
\displaystyle
=-\int_{\partial\Omega}\frac{\partial}{\partial\nu}\mbox{\boldmath $w$}^*_x\cdot\nabla G(z-y)\,dS(z)
\end{array}
$$
and the Alessandrini identity yields
$$\begin{array}{l}
\displaystyle
\,\,\,\,\,\,
\lim_{n\rightarrow\infty}\sum_{j}<(\Lambda_V-\Lambda_0)(v_n^j+H^j(\,\cdot\,,x)\vert_{\partial\Omega},v_n^j\vert_{\partial\Omega}>
\\
\\
\displaystyle
=\int_{D}V(z)(\mbox{\boldmath $w$}^*_x(z)+\nabla G(z-x)+\mbox{\boldmath $H$}(z,x))\cdot\nabla G(z-y)\,dz.
\end{array}
$$
Therefore, we obtain
$$\displaystyle
\lim_{n\rightarrow\infty}\sum_{j}\,
<(\Lambda_V-\Lambda_0)(v_n^j+H^j(\,\cdot\,,x))\vert_{\partial\Omega}, (v_n^j+H^j(\,\cdot\,,x))\vert_{\partial\Omega}>
=\nabla\cdot\mbox{\boldmath $w$}^*_x(x).
$$
This together with  (1.18) and (1.19) complets the proof of (1.23).

$\quad$

{\it\noindent Proof of {\rm(f)}.}
The formula (1.24) is derived from formulae (1.22), (1.23) and (1.6).

\noindent
$\Box$

\subsection{Uniqueness from local probing}

$\quad$

The lifting ${\mathcal I}(x,y)$ of ${\mathcal I}(x)$ satisfies the Laplace equation if one of two variables $x$ and $y$ is fixed, as shown (2.27) and (2.28).
Then by the weak unique continuation property for the solution of the Laplace equation, one immediately obtains

\begin{thm}
%\proclaim{\noindent Theorem 2.1.}
Let $U$ and $V$ be nonempty open subsets of $\Bbb R^3$ such that $U\cup V\subset\Omega$ and $\overline{U}\cap\overline{D}
=\overline{V}\cap\overline{D}=\emptyset$.
Then the values of ${\mathcal I}(x,y)$ for all $(x,y)\in(\Omega\setminus\overline{D})^2$  and thus
${\mathcal I}(x)={\mathcal I}(x,x)$ for all $x\in\Omega\setminus\overline{D}$ are uniquely determined by those of
${\mathcal I}(x,y)$ for all $(x,y)\in U\times V$.

%\endproclaim

\end{thm}

$\quad$

\noindent
The values of ${\mathcal I}(x,y)$ for all $(x,y)\in U\times V$
can be calculated from the following steps.

$\quad$

{\bf\noindent Step 1.}  Choose a point $x\in U$ and $\sigma\in N_x$ satisfying $\sigma(]0,\,1])\subset U$ (for example).

$\quad$

{\bf\noindent Step 2.}  Choose a needle sequence $\{v_n^j\}$ for $(x,\sigma)$.

$\quad$

{\bf\noindent Step 3.}  Observe $\Lambda_Vf_n^j$ for the input 
$f_n^j=v_n^j\vert_{\partial\Omega}$.

$\quad$

{\bf\noindent Step 4.}  Choose a point $y\in V$ and $\sigma'\in N_y$ satisfying $\sigma'(]0,\,1])\subset V$.

$\quad$

{\bf\noindent Step 5.}  Choose a needle sequence $\{v_n'^j\}$ for $(y,\sigma')$.

$\quad$

{\bf\noindent Step 6.}  Calculate ${\mathcal I}(x,y)$ via the formula (2.30):
$$\displaystyle
\lim_{n\rightarrow\infty}\sum_{j}\int_{\partial\Omega}(\Lambda_v-\Lambda_0)f_n^j\cdot v_n'^j(z)\,dS(z)={\mathcal I}(x,y).
$$

\noindent
Note that we use the observation data only in Step 3.  The other steps consist of calculation.

Finally we mention that the same statement as Theorem 2.1 for ${\mathcal I}(x,y)$ replaced with $\nabla\cdot\mbox{\boldmath $w$}_x(y)$ is also valid.
This is because of (2.29) and (trivial) harmonicity of $\nabla\cdot\mbox{\boldmath $w$}_x(y)$ with respect to $y\in\Omega\setminus\overline{D}$ for each fixed $x\in\Omega\setminus\overline{D}$.
The formula in Step 6 becomes (2.32):
$$\begin{array}{ll}
\displaystyle
\nabla\cdot\mbox{\boldmath $w$}_x(y)
&
\displaystyle
=-\lim_{n\rightarrow\infty}\sum_j
\int_{\partial\Omega}
(\Lambda_V-\Lambda_0)f_n^j\cdot \left(\frac{\partial G}{\partial z_j}(z-y)-v_n'^j(z)\right)\,dS(z).
\end{array}
$$

\section{Side B of IPS}

In this section, we consider the Side B of IPS.
The problem to be discussed here is as follows.

$\quad$

{\bf\noindent Problem 2.} Given $x\in\Omega$ and $\sigma\in N_x$, let $\{v_n^j\}$ be a needle sequence for $(x,\sigma)$.
When $\sigma\cap\overline{D}\not=\emptyset$, clarify the behaviour of the three sequences for each $j=1,2,3$ defined by
$$\displaystyle
{\mathcal I}_n(x,\sigma,\{v_n^j\})=
\left\{
\begin{array}{ll}
\displaystyle
<(\Lambda_V-\Lambda_0)v_n^j\vert_{\partial\Omega},v_n^j\vert_{\partial\Omega}>, & \text{Probe Method}
\\
\\
\displaystyle
-<(\Lambda_V-\Lambda_0)v_n^j\vert_{\partial\Omega}, G_n^j(\,\cdot\,,x)\vert_{\partial\Omega}>,& \text{Singular Sources Method}
\\
\\
\displaystyle
<(\Lambda_V-\Lambda_0)G_n^j(\,\cdot\,,x)\vert_{\partial\Omega},G_n^j(\,\cdot\,,x)\vert_{\partial\Omega}>. & \text{Completely Integrated Method}
\end{array}
\right.
$$

$\quad$

\noindent
We call the sequence $\{{\mathcal I}_n(x,\sigma,\{v_n^j\})\}$  for each method the {\it indicator sequence}.
Problem 2 is a question about the Side B of IPS.
Note that, in Problem 2 we do not take the summation of indicator sequence over the index $j=1,2,3$.  
The Side B of the Probe Method was firstly considered in \cite{INew} for an inverse obstacle problem governed by the Helmholtz equation.
However, there was no attempt to consider the Side B of the original Singular Sources Method \cite{P1}.  In \cite{IPS}, we considered the Side B
of the reformulated singular sources method for a prototype inverse obstacle problem governed by the Laplace equation.
However, its clarification was pending.  Recently as a byproduct of IPS developed in \cite{IPS2} and \cite{IPS3} we found
a simple solution about the Side B of the singular sources method under the IPS framework for impenetrable obstacle cases.
It is a reduction to the Side B of both the Probe Method and Completely Integrated Method.  This means, first of all, one has to solve Problem 2
for Probe Method together with Completely Integrated Method.

First we prepare the blowing up property of the needle sequence on the needle.

$\quad$

{\bf\noindent Proposition 3.1.}
%\proclaim{\noindent Proposition 3.1}   
{\it Given $x\in\Omega$ and $\sigma\in N_x$ let $\{v_n^j\}$ be an arbitrary needle sequence for $(x,\sigma)$.

\noindent
{\rm (a)}   Let $V$ be an arbitrary finite and open cone vertex at $x$.  We have
$$\displaystyle
\lim_{n\rightarrow\infty}\Vert v_n^j\Vert_{L^2(V\cap\Omega)}=\infty.
$$

\noindent
{\rm (b)}  Fix a point $z\in\sigma\setminus\{\sigma(0)\}$.  Let $B$ be an arbitrary open ball centered at $z$.
We have
$$\displaystyle
\lim_{n\rightarrow\infty}\Vert v_n^j\Vert_{L^2(B\cap\Omega)}=\infty.
$$
}

%\endproclaim

$\quad$

\noindent
Note that, in this proposition,  
we do not choose $\{v_n^j\}$ given by $v_n^j=\frac{\partial v_n}{\partial z_j}$, where $\{v_n\}$ is the original needle sequence in \cite{IProbe}.

{\it\noindent Proof.}
The proof of {\rm (a)} is a consequence of Lemma A of the Appendix, convergence property of $\{v_n^j\}$ in $\Omega\setminus\sigma$ 
and Fatou's lemma.  See Lemma 2.1 in \cite{INew}.

The proof of {\rm (b)} almost follows that of the corresponding fact to the original needle sequence for the Helmholtz equation, see Lemma 2.2 in \cite{INew}.
However, for readers convenience, we present here its argument since the Laplace equation case enables us to skip some of technical parts.
The proof is divided into the following five steps.

\noindent
(1). Choose a small open ball $B_1$ centered at $z$ in such a way that $B_1\subset B\cap\Omega$.

\noindent
(2). Choose a small open ball $B_2$ centered at $z$ and contained in $\overline{B_1}$ 
 in such a way that  $B_2\cap\sigma$ consists of a single connected piece of $\sigma$ passing through the point $z$.
 
 \noindent
 (3). Choose a nonempty open subset $S_0$ of $\partial B_2$ such that the intersection of $S_0$ and the piece of $\sigma$ connecting $z$ and $x$ consists
 of a singe point.

 \noindent
Then one gets
\begin{equation}
%$$
\displaystyle
\Vert v_n^j\Vert_{L^2(B\cap\Omega)}\ge \Vert v_n^j\Vert_{L^2(B_1)}\ge C_1\Vert v_n^j\Vert_{H^1(B_2)}\ge C_2\Vert v_n^j\Vert_{L^2(\partial B_2)}
\ge C_2\Vert v_n^j\Vert_{L^2(S_0)}.
\tag {3.1}
%$$
\end{equation}

\noindent
(4). Choose an open subset $U$ of $\Bbb R^3$ with a sufficiently smooth boundary such that $\overline{U}\subset\Omega$, $x\in U$,
$S_0\subset\partial U$ and $\partial U\setminus S_0\subset\Omega\setminus\sigma$.

\noindent
Then, by the convergence of $\{v_n^j\}$ in $\Omega\setminus\sigma$ 
\begin{equation}
%$$
\displaystyle
\sup_n\,\Vert v_n^j\Vert_{L^2(\partial U\setminus S_0)}=C_3<\infty.
\tag {3.2}
%$$
\end{equation}

\noindent
(5).  Choose an open ball $B_3$ centered at $x$ with $\overline{B_3}\subset U$.

\begin{equation}
%$$
\displaystyle
\Vert v_n^j\Vert_{L^2(B_3)}\le C_4\Vert v_n^j\Vert_{L^2(\partial U)}\le C_4(\Vert v_n^j\Vert_{L^2(\partial U\setminus S_0)}+\Vert v_n^j\Vert_{L^2(S_0)}).
\tag {3.3}
%$$
\end{equation}

\noindent
From (3.1), (3.2) and (3.3), we obtain
$$\displaystyle
\Vert v_n^j\Vert_{L^2(B\cap\Omega)}\ge C_2C_4^{-1}\Vert v_n^j\Vert_{L^2(B_3)}-C_2C_3.
$$
This together with {\rm(a)} yields {\rm (b)}.

\noindent
$\Box$

$\quad$

\noindent
Proposition 3.1 together with the convergence property of needle sequence implies: from the behaviour of needle sequence $\{v_n^j\}$ 
for $(x,\sigma)$ one can recover the full geometry of needle $\sigma$ by the formula:
$$\displaystyle
\sigma(]0,\,1])
=\{z\in\Omega\,\vert\, \text{$\lim_{n\rightarrow\infty}\Vert v_n^j\Vert_{L^2(B\cap \Omega)}=\infty$ for any open ball $B$ centered at $z$}\,\}.
$$

The following proposition is the key for Side B of IPS.

$\quad$

{\bf\noindent Proposition 3.2.}
%\proclaim{\noindent Proposition 3.2.}
{\it
Given $x\in\Omega$ and $\sigma\in N_x$ let $\{v_n^j\}$ be an arbitrary needle sequence for $(x,\sigma)$.  
Let $U$ be an arbitrary non empty open subset of $\Bbb R^3$ with $\overline{U}\subset\Omega$ with Lipschitz boundary satisfying one of 
two cases {\rm (i)} and {\rm (ii)} listed below:

\noindent
{\rm (i)} $x\in\overline{U}$;

\noindent
{\rm(ii)} $x\in\Omega\setminus\overline{U}$ and $\sigma\cap U\not=\emptyset$.

\noindent
Then, we have
\begin{equation}
%$$
\displaystyle
\lim_{n\rightarrow\infty}\Vert v_n^j\Vert_{L^2(U)}=\infty
\tag {3.4}
%$$
\end{equation}
and
\begin{equation}
%$$
\displaystyle
\lim_{n\rightarrow\infty}
\frac{\Vert v_n^j\Vert_{L^1(U)}}{\Vert v_n^j\Vert_{L^2(U)}}=0.
\tag {3.5}
%$$
\end{equation}

}

%\endproclaim

$\quad$

{\it\noindent Proof.}
First we prove (3.4).
If $x\in\overline{U}=\partial U\cup U$, one can take an nonempty open finite cone $V$ with the vertex at $x$ such that $V\subset U$.  
Then {\rm (a)} of Proposition 3.1 yields the desired conclusion.  If $x\in\Omega\setminus\overline{U}$ 
and $\sigma\cap U\not=\emptyset$, then choose a point $z\in \sigma\cap U$ and take an open ball $B$ centered at $z$ with $B\subset U$.
Then from {\rm (b)} of Proposition 3.1 one obtains the same conclusion.

The proof of (3.5) is as follows.
Given $\epsilon>0$ set
$$\displaystyle
\sigma_{\epsilon}=\{x\in\Bbb R^3\,\vert\,\text{dust}\,(x,\sigma)<\epsilon\}.
$$
We have
$$\begin{array}{ll}
\displaystyle
\Vert v_n^j\Vert_{L^1(U)}
&
\displaystyle
\le\Vert v_n^j\Vert_{L^1(U\cap\sigma_{\epsilon})}+\Vert v_n^j\Vert_{L^1(U\setminus\sigma_{\epsilon})}
\\
\\
\displaystyle
&
\displaystyle
\le \vert U\cap\sigma_{\epsilon}\vert^{\frac{1}{2}}\Vert v_n^j\Vert_{L^2(U)}+\vert U\vert^{\frac{1}{2}}\Vert v_n^j\Vert_{L^2(U\setminus\sigma_{\epsilon})}.
\end{array}
$$
Thus one gets
\begin{equation}
%$$
\displaystyle
\frac{\Vert v_n^j\Vert_{L^1(U)}}{\Vert v_n^j\Vert_{L^2(U)}}\le\vert U\cap\sigma_{\epsilon}\vert^{\frac{1}{2}}
+\vert U\vert^{\frac{1}{2}}\frac{\Vert v_n^j\Vert_{L^2(U\setminus\sigma_{\epsilon})}}{\Vert v_n^j\Vert_{L^2(U)}}.
\tag {3.6}
%$$
\end{equation}
Since $U\setminus\sigma_{\epsilon}\subset\Omega\setminus\sigma$, one has
$$\displaystyle
\lim_{n\rightarrow\infty}\Vert v_n^j\Vert_{L^2(U\setminus\sigma_{\epsilon})}=\left\Vert \frac{\partial G}{\partial z_j}(\,\cdot\,-x)\right\Vert_{L^2(U\setminus\sigma_{\epsilon})}<\infty.
$$
A combination of this, (3.4) and (3.6) yields
$$\displaystyle
\limsup_{n\rightarrow\infty}
\frac{\Vert v_n^j\Vert_{L^1(U)}}{\Vert v_n^j\Vert_{L^2(U)}}\le\vert U\cap\sigma_{\epsilon}\vert^{\frac{1}{2}}.
$$
Since $\vert U\cap\sigma_{\epsilon}\vert\rightarrow 0$ as $\epsilon\rightarrow 0$, we conclude the validity of  (3.5).

\noindent
$\Box$

$\quad$

Here we impose a jump condition on $V(z)$ which is stronger than the previous one in Theorems 1.1, 1.3 and 1.4.

$\quad$

{\bf\noindent Definition 3.1.}
We say the potential $V(z)$ has positive/negative jump on $D$ if there
exists a positive constant $C$ such that $\pm V(z)\ge C$ a. e. $z\in D$.
In what follows, unless otherwise specified, the phrase {\it on $D$} is omitted.

$\quad$

The Side B of IPS which is an answer to Problem 2 is as follows.

\begin{thm}
%\proclaim{\noindent Theorem 3.1.}
Let $x\in\Omega$ and $\sigma\in N_x$.  Assume that one of the two cases {\rm (i)} and {\rm (ii)} listed below is satisfied:

\noindent
{\rm (i)}  $x\in\overline{D}$;

\noindent
{\rm (ii)}  $x\in\Omega\setminus\overline{D}$ and $\sigma\cap D\not=\emptyset$.

\noindent
Then for any needle sequence $\{v_n^j\}$ for $(x,\sigma)$ we have
\begin{equation}
%$$
\displaystyle
\lim_{n\rightarrow\infty}\,{\mathcal I}_n(x,\sigma,\{v_n^j\}) =
\left\{
\begin{array}{ll}
\displaystyle
\infty 
&
\text{if $V(z)$ has positive jump,}
\\
\\
\displaystyle
-\infty
&
\displaystyle
\text{if $V(z)$ has negative jump.}
\end{array}
\right.
\tag {3.7}
%$$
\end{equation}

%\endproclaim
\end{thm}

{\it\noindent Proof.}  
We follow an approach presented in \cite{IR} for the Enclosure Method.
It means that,  we employ the Alessandrini identity, that is
$$\displaystyle
<(\Lambda_V-\Lambda_0)f,f>=\int_DV(z)(u-v)v\,dz+\int_D V(z)v^2\,dz.
$$
where $u$ is the solution of (1.1) with $f=v\vert_{\partial\Omega}$ and $v$ satisfies the Laplace equation in the whole domain $\Omega$.
Assume that $\pm V(z)\ge C$ a.e. $z\in D$ for a positive constant $C$.  We have
$$\displaystyle
\pm <(\Lambda_V-\Lambda_0)f,f>\ge C\Vert v\Vert_{L^2(D)}^2-\int_DV(z)(u-v)v\,dz.
$$
Here by Lemma 3.1 in \cite{IR}, we have
$$\displaystyle
\Vert u-v\Vert_{L^2(D)}\le C'\Vert v\Vert_{L^1(D)}.
$$
Thus one gets
\begin{equation}
%$$
\displaystyle
\pm <(\Lambda_V-\Lambda_0)f,f>\ge C\Vert v\Vert_{L^2(D)}^2-C''\Vert v\Vert_{L^2(D)}\Vert v\Vert_{L^1(D)}.
\tag {3.8}
%$$
\end{equation}

First consider the case when the indicator sequence has the form for the Probe Method:
$$
\displaystyle
{\mathcal I}_n(x,\sigma,\{v_n^j\})=<(\Lambda_V-\Lambda_0)v_n^j\vert_{\partial\Omega},v_n^j\vert_{\partial\Omega}>.
$$
By (3.8) with $f=v_n^j\vert_{\partial\Omega}$ and the assumption on $V$, it holds that
$$\displaystyle
\pm{\mathcal I}_n(x,\sigma,\{v_n^j\})\ge C\Vert v_n^j\Vert^2_{L^2(D)}\left(1-C'\frac{\Vert v_n^j\Vert_{L^1(D)}}{\Vert v_n^j\Vert_{L^2(D)}}\right).
$$
Then, applying Proposition 3.2 with $U=D$, we obtain the desired conclusion.

Next consider the case when the indicator sequence has the form for the Completely Integrated Method.
In this case
$$\displaystyle
{\mathcal I}_n(x,\sigma,\{v_n^j\})=<(\Lambda_V-\Lambda_0)G_n^j(\,\cdot\,,x)\vert_{\partial\Omega},G_n^j(\,\cdot\,,x)\vert_{\partial\Omega}>.
$$
By (1.15), this is another expression
$$\displaystyle
{\mathcal I}_n(x,\sigma,\{v_n^j\})=<(\Lambda_V-\Lambda_0)(v_n^j+H(\,\cdot\,,x))\vert_{\partial\Omega},(v_n^j+H(\,\cdot\,,x))\vert_{\partial\Omega}>.
$$
Then from (3.8) with $f=(v_n^j+H^j(\,\cdot\,,x))\vert_{\partial\Omega}$
one gets
$$\begin{array}{ll}
\displaystyle
\pm{\mathcal I}(x,\sigma,\{v_n^j\})
&
\displaystyle
\ge C\Vert v_n^j+H^j(\,\cdot\,x)\Vert^2_{L^2(D)}-C'\Vert v_n^j+H^j(\,\cdot\,x)\Vert_{L^2(D)}\Vert v_n^j+H^j(\,\cdot\,x)\Vert_{L^1(D)}.
\end{array}
$$
Here, by Proposition 3.2 with $U=D$, we have, as $n\rightarrow\infty$
$$\displaystyle
\Vert v_n^j+H^j(\,\cdot\,,x)\Vert_{L^2(D)}\ge \Vert v_n^j\Vert_{L^2(D)}-\Vert H^j(\,\cdot\,,x)\Vert_{L^2(D)}\rightarrow\infty
$$
and
$$\displaystyle
\frac{\Vert v_n^j+H^j(\,\cdot\,x)\Vert_{L^1(D)}}{\Vert v_n^j+H^j(\,\cdot\,x)\Vert_{L^2(D)}}
\le
\frac{\Vert v_n^j\Vert_{L^1(D)}+\Vert H^j(\,\cdot\,x)\Vert_{L^1(D)}}{\Vert v_n^j\Vert_{L^2(D)}-\Vert H^j(\,\cdot\,x)\Vert_{L^2(D)}}
\rightarrow 0.
$$
Hereafter we take the same course as above and obtain the conclusion.

Finally consider the indicator sequence for the Singular Sources Method.
The argument was taken from \cite{IPS2} and is purely algebraic.
It is a direct consequence of the formula
\begin{equation}
%$$
\begin{array}{l}
\displaystyle
\,\,\,\,\,\,
-<(\Lambda_V-\Lambda_0)v_n\vert_{\partial\Omega},G_n^j(\,\cdot\,,x)\vert_{\partial\Omega}>
\\
\\
\displaystyle
=\frac{1}{2}\left(
<(\Lambda_V-\Lambda_0)v_n^j\vert_{\partial\Omega}, v_n^j\vert_{\partial\Omega}>+<(\Lambda_V-\Lambda_0)G_n^j(\,\cdot\,,x)\vert_{\partial\Omega},G_n^j(\,\cdot\,,x)\vert_{\partial\Omega}>
\right.
\\
\\
\displaystyle
\,\,\,
\left.-<(\Lambda_V-\Lambda_0)\nabla G(\,\cdot\,-x)\vert_{\partial\Omega},\nabla G(\,\cdot\,-x)\vert_{\partial\Omega}>\right).
\end{array}
\tag {3.9}
%$$
\end{equation}

\noindent
$\Box$

$\quad$

\noindent
Note that using integration by parts, one gets also the decomposition
\begin{equation}
%$$
\displaystyle
<(\Lambda_V-\Lambda_0)f,f>=-\int_{\Omega}\vert\nabla (u-v)\vert^2\,dz-\int_DV(z)\vert u-v\vert^2\,dz+\int_DV(z)\vert v\vert^2\,dz,
\tag {3.10}
%$$
\end{equation}
where $u$ is the solution of (1.1) with $f=v\vert_{\partial\Omega}$ and $v$ satisfies the Laplace equation in the whole domain $\Omega$.
However, the second term of the right-hand side on (3.10) has minus sign, by this reason, however not all,  
it seems the Alessandrini identity is better than (3.10) for our purpose.

$\quad$

{\bf\noindent Remark 3.1.}
Consider the case when  $D$ has the decomposition $D=D_1\cup D_2$, where $\overline{D_1}\cap\overline{D_2}=\emptyset$ and each  $D_i$, $i=1,2$ is a nonempty
open subset of $\Bbb R^3$ with Lipschitz boundary; $V(z)$ has positive jump on $D_1$ and
negative jump on $D_2$.  Theorem 3.1 is still valid if one modifies the statement (3.7) as below:
$$\displaystyle
\lim_{n\rightarrow\infty}\,{\mathcal I}_n(x,\sigma,\{v_n^j\}) =
\left\{
\begin{array}{ll}
\displaystyle
\infty 
&
\text{if $V(z)$ has positive jump on $D_1$ and $\sigma\cap\overline{D_2}=\emptyset$,}
\\
\\
\displaystyle
-\infty
&
\displaystyle
\text{if $V(z)$ has negative jump on $D_2$ and $\sigma\cap\overline{D_1}=\emptyset$.}
\end{array}
\right.
$$
The proof is almost same.  This fact corresponds to Theorem 3.2 (and Theorem 3.3) of  \cite{IPS3} where both of 
the sound-soft and -hard obstacles are embedded in a homogeneous background medium.
Comparing both the statements we see that $D_1$ and $D_2$ correspond to being sound-soft and sound-hard, respectively
Note that we are considering $\Lambda_V-\Lambda_0$ not $\Lambda_0-\Lambda_V$.

$\quad$

{\bf\noindent Remark 3.2}. Let $x\in\Omega$ and let  $\sigma\in N_x$.
If $\sigma$ satisfis $\sigma\cap\overline{D}=\emptyset$, then for an arbitrary needle sequence $\{v_n^j\}$ for $(x,\sigma)$ with a fixed $j$
we have the convergence
$$\displaystyle
\lim_{n\rightarrow\infty}\,<(\Lambda_V-\Lambda_0)v_n^j\vert_{\partial\Omega},v_n^j\vert_{\partial\Omega}>
=\int_DV(z)(\frac{\partial G}{\partial z_j}(z-x)+(\mbox{\boldmath $w$}_x)_j(z))\frac{\partial G}{\partial z_j}(z-x)\,dz
$$
and
$$\displaystyle
\lim_{n\rightarrow\infty}\,<(\Lambda_V-\Lambda_0)G_n^j(\,\cdot\,,x),G_n^j(\,\cdot\,,x)\vert_{\partial\Omega}>
=\int_{D}V(z)(G^j(z,x)+(\mbox{\boldmath $w$}_x^*)_j(z))G^j(z,x)\,dz.
$$
The proof of these are the same as (1.3) and  (1.20) of Theorems 1.1 and 1.4, respectively.
Thus from (3.9) one gets also the limit
$$\displaystyle
-\infty<-\lim_{n\rightarrow\infty}\,<(\Lambda_V-\Lambda_0)v_n\vert_{\partial\Omega},G_n^j(\,\cdot\,,x)\vert_{\partial\Omega}><\infty.
$$

$\quad$

As a direct corollary of Theorem 3 together with Remark 3.2
we obtain a characterization of obstacle in terms of the blowing up property of indicator sequences with {\it fixed} $j$.
In what follows, we assume that $\Omega\setminus\overline{D}$ is connected.  
This ensures that for all $x\in\partial D$ there exists a needle
$\sigma\in N_x$ such that $\sigma(t)\in\Omega\setminus\overline{D}$ for all $t\in]0,\,1[$.

$\quad$

{\bf\noindent Corollary 3.1.}
%\proclaim{\noindent Corollary 3.1.} 
{\it
A point $x\in\Omega$ belongs to $\overline{D}$ if and only if
for arbitrary $\sigma\in N_x$ and needle sequence $\{v_n^j\}$ for $(x,\sigma)$
$$\displaystyle
\lim_{n\rightarrow\infty}\,{\mathcal I}_n(x,\sigma,\{v_n^j\})=
\left\{\begin{array}{ll}
\displaystyle
\infty
&
\displaystyle
\text{if $V(z)$ has positive jump,}
\\
\\
\displaystyle
-\infty
&
\displaystyle
\text{if $V(z)$ has negative jump.}
\end{array}
\right.
$$
}

%\endproclaim

$\quad$

\noindent
This result includes a criterion for a target distinction by the signature of the limit.
And this characterization of obstacle itself by using the {\it divergence} of a sequence calculated from the Dirichlet-to-Neumann map,
has a contrast in idea to Kirsch's factorization method \cite{Ki} developed for inverse medium scattering problem since it is the one by using the {\it convergence} of an infinite series calculated from the far field operator.

From Corollary 3.1 one gets also a simpler characterization:

$\quad$

{\bf\noindent Corollary 3.2.}
%\proclaim{\noindent Corollary 3.2.}
{\it
Under the assumption that $V(z)$ has positive or negative jump, a point $x\in\Omega$ belongs to $\overline{D}$ if and only if
for arbitrary $\sigma\in N_x$ and needle sequence $\{v_n^j\}$ for $(x,\sigma)$
$$\displaystyle
\lim_{n\rightarrow\infty}\vert{\mathcal I}_n(x,\sigma,\{v_n^j\})\vert=\infty.
$$
}

%\endproclaim

$\quad$

This is an obstacle characterization by using Side B of IPS.  However, this characterization requires us to check all the needles and needle sequences.
To overcome this one can give an alternative, however, equivalent version of this corollary.  That is

$\quad$

{\bf\noindent Corollary 3.3.}
%\proclaim{\noindent Corollary 3.3.}   
{\it
Under the same assumption on $V(z)$ as Corollary 3.2, it holds that:
a point $x\in\Omega$ belongs to $\Omega\setminus\overline{D}$ if and only if
there exists a needle $\sigma\in N_x$ and for some needle sequence $\{v_n^j\}$ for $(x,\sigma)$
$$\displaystyle
\liminf_{n\rightarrow\infty}\vert{\mathcal I}_n(x,\sigma,\{v_n^j\})\vert<\infty.
$$
}

%\endproclaim

$\quad$

However, clearly, we can say more that

$\quad$

{\bf\noindent Corollary 3.4.}
%\proclaim{\noindent Corollary 3.4.} 
{\it
Under the same assumption on $V(z)$ as {\rm Corollary 3.2}, it holds that: a point $x\in\Omega$ belongs to $\Omega\setminus\overline{D}$ if and only if
there exists a needle $\sigma\in N_x$ and for some needle sequence $\{v_n^j\}$ for $(x,\sigma)$ the limit
$$\displaystyle
-\infty<\lim_{n\rightarrow\infty}{\mathcal I}_n(x,\sigma,\{v_n^j\})<\infty
$$
exists.
}

%\endproclaim

$\quad$

\noindent
Needles to say, this criterion requires us a decision whether given sequence is convergent or not.

Corollary 3.4 gives a characterization {\it outside} obstacle by using the {\it convergence} of a sequence calculated from the Dirichlet-to-Neumann map,
which has a contrast in idea to Kirsch's factorization method \cite{Ki} since it is the one by using the {\it divegence} of an infinite series calculated from the far field operator.

$\quad$

{\bf\noindent Remark 3.3.}
Since $x\in\Omega\setminus\overline{D}$ is independent of $j$, 
under the same assumption on $V(z)$ as {\rm Corollary 3.2}, we conclude also that the following statements for $j=1,2,3$ are equivalent to each other:
there exists a needle $\sigma\in N_x$ and for some needle sequence $\{v_n^j\}$ for $(x,\sigma)$ the limit
$$\displaystyle
-\infty<\lim_{n\rightarrow\infty}{\mathcal I}_n(x,\sigma,\{v_n^j\})<\infty
$$
exists.

\section{Remarks on Future Directions of IPS}

In this paper, we derived a singular sources method by utilizing the probe method for the Schr\"odinger equation and identifying the corresponding IPS function.
The derivation of the IPS function is based on Green's representation.  There would be no difficulty in covering the case where the potential takes a nonzero constant value
outside the obstacle.  In particular, this is the case when  the Schr\"odinger equation in (1.1) is replaced with the inhomogeneous Helmholtz equation
$$\begin{array}{ll}
\displaystyle
\Delta u+k^2(1+V(x))u=0, & x\in\Omega,
\end{array}
$$ 
where $V$ satisfies (1.2) and $k>0$.
It would be interesting to apply IPS
to cases where the potential modeling the background medium is neither a constant nor real-valued function.

In \cite{IProbeElastic} (see also Subsection 2.2.3 of the review paper \cite{IReview}),
the probe method based on the needle sequences generated by two types of singular solutions for the Navier system has been
established.  One comes from the Kelvin matrix; the other, however, is not of that type.   Thus, it would be interesting to develop IPS for the Navier system.
This would involve reconsidering the probe method developed therein from an IPS perspective
and introducing the corresponding singular sources method within the IPS framework.

\section{Appendix}

\subsection{Divergence of an integral}

$\quad$

{\bf\noindent Lemma A.}
%\proclaim{\noindent Lemma A.}  
{\it
Let $\mbox{\boldmath $a$}$ be a constant unit vector and $V$ be an open and finite cone around vector $\mbox{\boldmath $a$}$ with vertex at the origin of the coordinate, that is
$$\displaystyle
V=V(\mbox{\boldmath $a$},\Theta)=\left\{z\in\Bbb R^3\,\vert\, z\cdot\mbox{\boldmath $a$}>\vert z\vert\cos\Theta\,\right\}\cap B
$$
where $\Theta\in\,]0,\frac{\pi}{2}[$ and $B$ is an open ball centered the origin.
Then, it holds that,  for an arbitrary unit vector $\mbox{\boldmath $b$}$
\begin{equation}
%$$
\displaystyle
\int_V\vert\nabla G(z)\cdot\mbox{\boldmath $b$}\vert^2\,dz=\infty.
\tag {{\rm A.1}}
%$$
\end{equation}
}
%\endproclaim

$\quad$

{\it\noindent Proof.}
First consider the case when $\mbox{\boldmath $a$}\cdot\mbox{\boldmath $b$}\not=0$.
Let $n=1,2,\cdots$ and
$$\displaystyle
V_n=\left\{z\in\Bbb R^3\,\vert\, z\cdot\mbox{\boldmath $a$}>\vert z\vert\cos(\frac{\Theta}{n})\right\}\cap B.
$$
We have $V_n\subset V$ and thus
\begin{equation}
%$$
\displaystyle
\int_V\vert\nabla G(z)\cdot\mbox{\boldmath $b$}\vert^2\,dz\ge\int_{V_n}\vert\nabla G(z)\cdot\mbox{\boldmath $b$}\vert^2\,dz=\lim_{\epsilon\rightarrow 0}\int_{\epsilon}^R\frac{dr}{r^2}\int_{S_n(\mbox{\boldmath $a$})}\vert\mbox{\boldmath $\omega$}\cdot\mbox{\boldmath $b$}\vert^2\,dS(\mbox{\boldmath $\omega$}),
\tag {\rm {A.2}}
%$$
\end{equation}
where $R$ is the radius of $B$ and $S_n(\mbox{\boldmath $a$})$ is a non empty open subset of $S^2$ defined by
$$\displaystyle
S_n(\mbox{\boldmath $a$})=\left\{\mbox{\boldmath $\omega$}\in S^2\,\vert\, \mbox{\boldmath$\omega$}\cdot\mbox{\boldmath $a$}>\cos(\frac{\Theta}{n})\,\right\}.
$$
We show that there exists a $n_0$ such that, for all $n\ge n_0$
\begin{equation}
%$$
\displaystyle
\inf_{\mbox{\boldmath $\omega$}\in S_n(\mbox{\boldmath $a$})}\vert\mbox{\boldmath $\omega$}\cdot\mbox{\boldmath $b$}\vert>0.
\tag {{\rm A.3}}
%$$
\end{equation}
If not so is, there exists a sequence $\{\mbox{\boldmath $\omega$}_n\}$ such that, each $\mbox{\boldmath $\omega$}_n\in S_{n}(\mbox{\boldmath $a$})$ and satisfies
$$\displaystyle
\vert\mbox{\boldmath $\omega$}_n\cdot\mbox{\boldmath $b$}\vert<\frac{1}{n}.
$$
Since $\vert\mbox{\boldmath $\omega$}_n-\mbox{\boldmath $a$}\vert^2<2(1-\cos\frac{\Theta}{n})\rightarrow 0$, we conclude $\mbox{\boldmath $a$}\cdot\mbox{\boldmath $b$}=0$.  This is a contradiction.
Therefore we have (A.3) and thus (A.2) yields (A.1).

Next consider the case when $\mbox{\boldmath $b$}\cdot\mbox{\boldmath $a$}=0$.  In this case,
one can choose another unit vector $\mbox{\boldmath $c$}$ and $\theta\in\,]0,\Theta[$ in such a way that $V(\mbox{\boldmath $c$},\theta)\subset V(\mbox{\boldmath $a$},\Theta)$ and $\mbox{\boldmath $b$}\cdot\mbox{\boldmath $c$}\not=0$.  Then we have
$$\displaystyle
\int_V\vert\nabla G(z)\cdot\mbox{\boldmath $b$}\vert^2\,dz\ge \int_{V(\mbox{\boldmath $c$},\theta)}\vert\nabla G(z)\cdot\mbox{\boldmath $b$}\vert^2\,dz=\infty.
$$

\noindent
$\Box$

\subsection{Blowing-up of the IPS function on $\partial\Omega$}

$\quad$

In this subsection we give a proof of

$\quad$

{\bf\noindent Proposition B.}
%\proclaim{\noindent Proposition B.}
{\it
Given $b\in\partial\Omega$, we have
$$\displaystyle
\lim_{x\rightarrow b}\,{\mathcal I}^1(x)=-\infty.
$$
}
%\endproclaim

$\quad$

Once we have this, then the following statement is deduced from {\rm(c)} of Theorems 1.1 and 1.2:
it holds that, for each $b\in\partial\Omega$
$$\displaystyle
\lim_{x\rightarrow b}\nabla\cdot\mbox{\boldmath $W$}_x(x)=-\infty.
$$

Proof of Proposition B is based on a crucial estimate.

$\quad$

{\bf\noindent Lemma B.1.}
%\proclaim{\noindent Lemma B.1.}
{\it
Let $W$ be a bounded domain with $C^{1,1}$-boundary of $\Bbb R^3$.  
Let $V\in L^{\infty}(W)$.  Let $v\in H^1(W)$ satisfy
$-\Delta v+V(z)v=0$ in $W$.  Assume that, $0$ is not a Dirichlet eigenvalue of $-\Delta+V(z)$ in $W$.
Then, there exists a positive constant $C$ independent of $v$ such that
\begin{equation}
%$$
\displaystyle
\Vert v\Vert_{L^2(W)}^2\le C\Vert v\vert_{\partial W}\Vert_{L^2(\partial W)}^2.
\tag {{\rm B.1}}
%$$
\end{equation}
}
%\endproclaim

$\quad$

\noindent
For a proof, apply the argument for the proof of (A.2) of Proposition A.1 in \cite{INew}.

$\quad$

{\it\noindent Proof of Proposition B.}
Applying (B.1) to the solution of (1.7), one gets
$$\displaystyle
\Vert\mbox{\boldmath $w$}_x^1\Vert_{L^2(\Omega)}^2\le C
\Vert\nabla G(\,\cdot\,-x)\vert_{\partial\Omega}\Vert_{L^2(\partial\Omega)}^2.
$$
Then, by (1.10) we obtain
$$\displaystyle
-{\mathcal I}^1(x)\ge \Vert\nabla^2 G(\,\cdot\,-x)\Vert_{L^2(\Bbb R^3\setminus\overline{\Omega})}^2-C'
\Vert\nabla G(\,\cdot\,-x)\vert_{\partial\Omega}\Vert_{L^2(\partial\Omega)}^2,
$$
where $C'$ is a positive constant.
Thus it suffices to prove that
\begin{equation}
%$$
\displaystyle
\lim_{x\rightarrow b}
\Vert\nabla^2 G(\,\cdot\,-x)\Vert_{L^2(\Bbb R^3\setminus\overline{\Omega})}^2
=\infty
\tag {{\rm B.2}}
%$$
\end{equation}
and
\begin{equation}
%$$
\displaystyle
\lim_{x\rightarrow b}
\frac{\Vert\nabla G(\,\cdot\,-x)\vert_{\partial\Omega}\Vert_{L^2(\partial\Omega)}^2}
{
\Vert\nabla^2 G(\,\cdot\,-x)\Vert_{L^2(\Bbb R^3\setminus\overline{\Omega})}^2
}=0.
\tag {{\rm B.3}}
%$$
\end{equation}
Ignoring a constant multiplication, we have
$$\displaystyle
\Vert \nabla^2 G(\,\cdot\,-x)\Vert_{L^2(\Bbb R^3\setminus\overline{\Omega})}^2
\sim
\int_{\Bbb R^3\setminus\overline{\Omega}}\frac{dz}{\vert z-x\vert^6}.
$$
Then from Lemma B.2 mentioned below we obtain (B.2) and (B.3).

\noindent
$\Box$

$\quad$

{\bf\noindent Lemma B.2.}
%\proclaim{\noindent Lemma B.2.} 
{\it  Assume that $\partial\Omega$ be Lipschitz.
We have, as $x\rightarrow b$ in $\Omega$
\begin{equation}
%$$
\displaystyle
\int_{\Bbb R^3\setminus\overline{\Omega}}\frac{dz}{\vert z-x\vert^6}
\ge \frac{C_1}{\vert x-b\vert^3}
\tag {{\rm B.4}}
%$$
\end{equation}
and
\begin{equation}
%$$
\displaystyle
\frac{\Vert\nabla G(\,\cdot\,-x)\vert_{\partial\Omega}\Vert_{L^2(\partial\Omega)}^2}
{
\Vert\nabla^2 G(\,\cdot\,-x)\Vert_{L^2(\Bbb R^3\setminus\overline{\Omega})}^2
}
\le C_2\sqrt{\vert x-b\vert},
\tag {{\rm B.5}}
%$$
\end{equation}
where $C_1$ and $C_2$ are positive constants.
}
%\endproclaim

$\quad$

{\it\noindent Proof.}  We follow the proof of Theorem C in \cite{IProbe}.
Choose the local coordinates at the point $b$ in such a way that
$$\displaystyle
(\Bbb R^3\setminus\overline{\Omega})\cap B(b,\delta)
=
\left\{z=b+A\left(\begin{array}{l}
\displaystyle
s
\\
\\
\displaystyle
t
\end{array}
\right)\,\left\vert\right.
\,
\varphi(s)<t\,,\vert s\vert^2+t^2<\delta^2
\right\},
$$
where $\delta=\delta(b)$ is a small positive constant, $(s,t)\in\Bbb R^2\times\Bbb R$, $\varphi\in C^{0,1}(\Bbb R^2)$ with compact support and satisfies $\varphi(0)=0$, and $A$ is a $3\times 3$-orthogonal matrix.

First of all we have
\begin{equation}
%$$
\begin{array}{ll}
\displaystyle
\int_{\Bbb R^3\setminus\overline{\Omega}}\frac{dz}{\vert z-x\vert^6}
&
\displaystyle
\ge
\int_{(\Bbb R^3\setminus\overline{\Omega})\cap B(b,\delta_b)}\frac{dz}{\vert z-x\vert^6}.
\end{array}
\tag {{\rm B.6}}
%$$
\end{equation}
Since
$$\displaystyle
\vert z-x\vert\le(\vert s\vert^2+t^2)^{\frac{1}{2}}
+\vert x-b\vert,
$$
one gets
\begin{equation}
%$$
\begin{array}{ll}
\displaystyle
\int_{(\Bbb R^3\setminus\overline{\Omega})\cap B(b,\delta_b)}\frac{dz}{\vert z-x\vert^6}
&
\displaystyle
\ge
C
\int_{\varphi(s)<t,\,\vert s\vert^2+t^2<\delta^2}\frac{ds dt}{\{(\vert s\vert^2+t^2)^{\frac{1}{2}}+\vert x-b\vert\}^6}
\\
\\
\displaystyle
&
\displaystyle
=\frac{C}{\vert x-b\vert^3}
\int_{\varphi(\vert x-b\vert s)<\vert x-b\vert t,\,\vert s\vert^2+t^2<\delta^2\vert x-b\vert^{-2}}\frac{ds dt}{\{(\vert s\vert^2+t^2)^{\frac{1}{2}}+1\}^6}.
\end{array}
\tag {{\rm B.7}}
%$$
\end{equation}
Here we use the change of variables
$$\begin{array}{ll}
\displaystyle
\left(\begin{array}{l}
\displaystyle
s
\\
\\
\displaystyle
t
\end{array}
\right)
=\rho\,\mbox{\boldmath $\omega$},
&
\rho>0,\,\mbox{\boldmath $\omega$}=(\omega',\omega_3)\in S^2, \omega'=(\omega_1,\omega_2).
\end{array}
$$
One gets
\begin{equation}
%$$
\begin{array}{l}
\,\,\,\,\,\,
\displaystyle
\int_{\varphi(\vert x-b\vert s)<\vert x-b\vert t,\,\vert s\vert^2+t^2<\delta^2\vert x-b\vert^{-2}}\frac{ds dt}{\{(\vert s\vert^2+t^2)^{\frac{1}{2}}+1\}^6}
\\
\\
\displaystyle
=\int_{S^2}\int_{\varphi(\vert x-b\vert \rho\omega')<\vert x-b\vert\rho\omega_3,\,\rho^2<\delta^2\vert x-b\vert^{-2}}\frac{\rho^2d\rho d\mbox{\boldmath $\omega$}}{(\rho+1)^6}.
\end{array}
\tag {{\rm B.8}}
%$$
\end{equation}
Choose a positive constant $K$ such that $\vert \varphi(s)\vert\le K\vert s\vert$.
Let $\mbox{\boldmath $\omega$}$ satisfies $\omega_3>K\vert\omega'\vert$.  Then we have
$$\begin{array}{ll}
\displaystyle
\varphi(\vert x-b\vert\rho\omega')
&
\displaystyle
\le K\vert x-b\vert\rho\vert\omega'\vert
\\
\\
\displaystyle
&
\displaystyle
<\vert x-b\vert\rho\omega_3.
\end{array}
$$
Thus one gets
\begin{equation}
%$$
\begin{array}{l}
\displaystyle
\int_{S^2}\int_{\varphi(\vert x-b\vert \rho\omega')<\vert x-b\vert\rho\omega_3,\,\rho^2<\delta^2\vert x-b\vert^{-2}}\frac{\rho^2d\rho\,d\mbox{\boldmath $\omega$}}{(\rho+1)^6}
\\
\\
\displaystyle
\ge
\int_{\omega\in S^2,\,\omega_3>K\vert\omega'\vert}d\mbox{\boldmath $\omega$}
\int_0^{\frac{\delta}{\vert x-b\vert}}
\frac{\rho^2}{(\rho+1)^6}d\rho.
\end{array}
\tag {{\rm B.9}}
%$$
\end{equation}
Since
$$\displaystyle
\int_{\omega\in S^2,\,\omega_3>K\vert \omega'\vert}d\mbox{\boldmath $\omega$}>0
$$
and
$$\displaystyle
0<\int_0^{\infty}\frac{\rho^2}{(\rho+1)^6}d\rho<\infty,
$$
one obtains, as $x\rightarrow b$
$$\displaystyle
\int_{\omega\in S^2,\,\omega_3>K\vert\omega'\vert}d\mbox{\boldmath $\omega$}
\int_0^{\frac{\delta}{\vert x-b\vert}}
\frac{\rho^2}{(\rho+1)^6}d\rho
\ge C',
$$
where $C'$ is a positive constant.
Therefore, from this together with (B.6) to (B.9), we obtain (B.4).

The proof of (B.5) is as follows.
By Theorem 1.5.1.10 in \cite{Gr}, taking a large open ball $B$,  we have, for all $\epsilon\in\,]0,\,1[$
$$\displaystyle
\Vert \nabla G(\,\cdot\,-x)\Vert_{L^2(\partial\Omega)}^2\le K(\epsilon\Vert\nabla^2G(\,\cdot\,-x)\Vert_{L^2(B\setminus\overline\Omega)}^2+
\epsilon^{-1}\Vert\nabla G(\,\cdot\,-x)\Vert_{L^2(B\setminus\overline\Omega)}^2),
$$
where $K$ is a positive constant independent of $\epsilon$ and $x\in\Omega$.
Thus one gets
\begin{equation}
%$$
\displaystyle
\frac{\Vert\nabla G(\,\cdot\,-x)\vert_{\partial\Omega}\Vert_{L^2(\partial\Omega)}^2}
{
\Vert\nabla^2 G(\,\cdot\,-x)\Vert_{L^2(\Bbb R^3\setminus\overline{\Omega})}^2
}
\le
K\left(\epsilon+\epsilon^{-1}
\frac{
\Vert\nabla G(\,\cdot\,-x)\Vert_{L^2(B\setminus\overline{\Omega)}}^2
}
{
\Vert\nabla^2 G(\,\cdot\,-x)\Vert_{L^2(\Bbb R^3\setminus\overline{\Omega})}^2
}
\right).
\tag {\rm{B.10}}
%$$
\end{equation}
Here we have, with the help of H\"older's inequality
$$\begin{array}{ll}
\displaystyle
\Vert\nabla G(\,\cdot\,-x)\Vert_{L^2(B\setminus\overline{\Omega)}}^2
&
\displaystyle
=C_1\int_{B\setminus\overline{\Omega}}\frac{dz}{\vert z-x\vert^4}
\\
\\
\displaystyle
&
\displaystyle
\le C_2\left(\int_{B\setminus\overline{\Omega}}
\frac{dz}{\vert z-x\vert^6}\right)^{\frac{2}{3}}
\\
\\
\displaystyle
&
\displaystyle
=C_3 \left(
\Vert\nabla^2 G(\,\cdot\,-x)\Vert_{L^2(B\setminus\overline{\Omega})}^2
\right)^{\frac{2}{3}}
\\
\\
\displaystyle
&
\displaystyle
=C_3\Vert\nabla^2 G(\,\cdot\,-x)\Vert_{L^2(B\setminus\overline{\Omega})}^{\frac{4}{3}},
\end{array}
$$
where $C_i$, $i=1,2,3$ are positive constants.
This yields
$$\displaystyle
\frac{
\Vert\nabla G(\,\cdot\,-x)\Vert_{L^2(B\setminus\overline{\Omega)}}^2
}
{
\Vert\nabla^2 G(\,\cdot\,-x)\Vert_{L^2(\Bbb R^3\setminus\overline{\Omega})}^2
}
\le\frac{C}{\Vert\nabla^2 G(\,\cdot\,-x)\Vert_{L^2(\Bbb R^3\setminus\overline{\Omega})}^{\frac{2}{3}}}.
$$
From this together with (B.4) one has, as $x\rightarrow b$
$$\displaystyle
\frac{
\Vert\nabla G(\,\cdot\,-x)\Vert_{L^2(B\setminus\overline{\Omega)}}^2
}
{
\Vert\nabla^2 G(\,\cdot\,-x)\Vert_{L^2(\Bbb R^3\setminus\overline{\Omega})}^2
}
\le C\vert x-b\vert.
$$
Thus from (B.10) one gets
$$\displaystyle
\frac{\Vert\nabla G(\,\cdot\,-x)\vert_{\partial\Omega}\Vert_{L^2(\partial\Omega)}^2}
{
\Vert\nabla^2 G(\,\cdot\,-x)\Vert_{L^2(\Bbb R^3\setminus\overline{\Omega})}^2
}
\le
K(\epsilon+\epsilon^{-1}C\vert x-b\vert).
$$
Thus, as $x\rightarrow b$, choosing a small $\epsilon$ in such a way that
$$\displaystyle
\epsilon=\epsilon^{-1}C\vert x-b\vert,
$$
we obtain (B.5).

\noindent
$\Box$

\noindent
Note that the proof of (B.5) presented above never makes use of any {\it upper bound} of $L^2$-norm of $\nabla G(\,\cdot\,-x)$ over $\partial\Omega$.

$\quad$

\end{document}